\documentstyle[12pt,a4]{amsart}

\numberwithin{equation}{section}

\newtheorem{thm}{Theorem}
\newtheorem{cor}{Corollary}
\newtheorem{lem}{Lemma}
\newtheorem{prop}{Proposition}

\begin{document}
\title[\null]{
A relation between the zeros of different two $L$-functions which have the Euler product
and functional equation}
\maketitle

\maketitle
\begin{center}
Masatoshi Suzuki
\smallskip

(December 12, 2004)
\end{center}

\begin{abstract}
As automorphic $L$-functions or Artin $L$-functions, 
several classes of $L$-functions have Euler products and functional equations. 
In this paper we study the zeros of $L$-functions 
which have the Euler products and functional equations. 
We show that 
there exists some relation between the zeros of the Riemann zeta-function 
and the zeros of such $L$-functions. 
As a special case of our results, we find the relations 
between the zeros of the Riemann zeta-function and 
the zeros of automorphic $L$-functions attached to elliptic modular forms 
or the zeros of Rankin-Selberg $L$-functions attached to two elliptic modular forms. 
\end{abstract}

%
%

\section{Introduction} 
Since the epoch making paper of G.B.Riemann ~\cite{Ri59}, 
the study of the zeros of the Riemann zeta-function and the other zeta-functions are  
one of the major fields of number theory. 
He introduced the analytic method to the theory of prime distribution, 
and clarified the relation between the primes and the zeros of the Riemann zeta-function.

After Riemann's paper, 
the accumulation of studies about the zeros of zeta-functions are enormous ones 
and many papers about them have been published every years. 
However almost of them have dealt with the zeros of single zeta-function, 
even if there exist some studies which deal with the zeros of a family of zeta-functions. 
An example is the studies related with the GUE-conjecture ~\cite{KaSa99}. 
It has given a new point of view to the study of the zeros of zeta-functions. 
This example tells us that 
studying the relation between the zeros of several different zeta-functions  
would give a new insight to the theory of the zeros of zeta-functions.    
Hence, in the present paper, we mainly concern to the relation between the zeros of different zeta-functions 
rather than individual properties of the zeros of single zeta-function. 
It is a continuation of the studies in ~\cite{Su04}.      

In ~\cite{Su04} the author showed that there exists some relation 
between the zeros of a $L$-function $L(s)$ belonging to the Selberg class 
and the zeros of an associated $L$-function $L_\chi(s)$ 
twisted by a primitive Dirichlet character $\chi$ 
which is a generalization of Linnik's result in~\cite{Li47}. 
Linnik's result is the asymptotic relation
\begin{equation}\label{101}
\sum_{{L(\rho,\chi)=0}\atop{0 < {\rm Re}(\rho) < 1}}
 \Gamma(\rho) x^{-\rho}
=\frac{1}{\tau(\overline{\chi})}
 \sum_{a=1}^q \overline{\chi}(a)
\sum_{{\zeta(\rho)=0}\atop{0 < {\rm Re}(\rho) < 1}}
 \Gamma(\rho) \left(x -2\pi i \frac{a}{q}\right)^{-\rho}
 +O(\log^2 x) 
\end{equation}
as $x \to +0$, where $\zeta(s)$ is the Riemann zeta function and 
$L(s,\chi)$ is the Dirichlet $L$-function attached to a primitive Dirichlet character $\chi$ mod $q$. 
The relation $\eqref{101}$ suggests that 
there exists a relation between the zeros of $\zeta(s)$ and $L(s,\chi)$.  
In fact Sprindzuk showed that, under the original Riemann hypothesis (RH), 
some properties of the zeros of $\zeta(s)$ are equivalent to 
the RH for Dirichlet $L$-functions by using $\eqref{101}$ 
in ~\cite{Sp75}. 
The author generalized the Sprindzuk type result to the case of the above pair  $L(s)$ and $L_\chi(s)$ in ~\cite{Su04}. 
Further studies on Sprinduzuk's work in ~\cite{Sp75}, see Fujii ~\cite{Fu88,Fu92}. 
The aim of this paper is to generalize the relation $\eqref{101}$ to the class of $L$-functions as wide as possible. 
The author believes that such results will give a new view point to the theory of $L$-functions.   

Our starting point is the observation about the zeros of principal automorphic $L$-functions. 
We refer to ~\cite{RuSa96} for notations and properties of them. 
Let $\pi=\otimes_p \pi_p$ be an irreducible cuspidal automorphic representation of $GL_N({\Bbb A}_{\Bbb Q})$ 
with unitary central character. 
The associated $L$-function $L(s,\pi)$ is given by a product of local factors $L(s,\pi_p)$.  
Except for a finite set of primes, $\pi_p$ is unramified. 
The local factors $L(s,\pi_p)$ for unramified primes are given by 
\begin{equation}\label{102}
L(s,\pi_p)=\prod_{j=1}^N (1-\alpha_{\pi}(p,j)p^{-s})^{-1}
\end{equation}
where $\alpha_{\pi}(p,j)$ are the eigen values of 
the semi-simple conjugacy class $\{A_{\pi}(p)\} \in GL_N(\Bbb C)$ associated to $\pi_p$. 
The generalized Ramanujan conjecture for cuspidal automorphic representation $\pi$ 
asserts that $|\alpha_{\pi}(p,j)|=1$ for unramified $p$. 
\smallskip

\noindent
{\bf Observation 1.} Under the general Ramanujan conjecture, 
the set of all zeros of $L(s,\pi_p)^{-1}$ for the unramified prime $p$  
is a union of $n$-piece translations of the zeros of $\zeta_{p}^{-1}(s)=1-p^{-s}$. 
\smallskip

The logarithmic derivative of $L(s,\pi)$ is written as  
\begin{equation}\label{103}
-\frac{L^\prime}{L}(s,\pi)=\sum_{n=1}^\infty \frac{\Lambda_\pi(n)}{n^s}
\end{equation}
where $\Lambda_\pi(n)=\Lambda(n)a_\pi(n)$, $\Lambda(n)=\log p$ if $n=p^m$ and zero otherwise, and  
\begin{equation}\label{104}
a_\pi(p^m)=\sum_{j=1}^N \alpha_{\pi}(p,j)^m. 
\end{equation}

By an ``explicit formula''  we usually mean an equation that 
represents the information of the Euler product and the functional equation 
in terms of an explicit relation between the zeros of $L(s,\pi)$ and $\Lambda_\pi(p^m)$. 
Let $h \in C_0^\infty({\Bbb R}_+)$ be a smooth compactly supported function, 
and let $\widehat{h}(s)=\int_0^\infty h(u)u^{s-1} du$ be the Mellin transform of $h$. 
Then
\begin{equation}\label{105}
\delta_\pi  \, \widehat{h}(0) -\sum_{L^\ast(\rho,\pi)=0}\widehat{h}(\rho) + \delta_\pi \, \widehat{h}(1) 
=\sum_{p} W_\pi(h,p) + W_{\pi}(h,\infty),
\end{equation}
where $\delta_\pi=1$ if $\pi$ corresponds to $\zeta(s)$, and zero otherwise, 
\begin{equation}\label{106}
W_\pi(h,p)=\sum_{m=1}^\infty (\Lambda_\pi(p^m) h(p^m) +\overline{\Lambda_\pi(p^m)} p^{-m}h(p^{-m})). 
\end{equation}
On the other hand, Poisson's summation formula yields   
\begin{equation}\label{107}
W_\pi(h,p)=\sum_{L(\rho,\pi_p)^{-1}=0} ( \widehat{h}(\rho) + \widehat{h}(1-\overline{\rho})).
\end{equation}

\noindent
{\bf Observation 2.} 
Combining $\eqref{105}$ with $\eqref{107}$, 
we obtain a relation between the zeros of $L(s,\pi)$ and the zeros of $L(s,\pi_p)^{-1}$. 
(We deal with similar things more precisely in $\S 8.3$.)
\smallskip

From {\bf Observation 1} and {\bf Observation 2}, 
we can see a possibility of the generalization of $\eqref{101}$ 
to principal automorphic $L$-functions. 
Further we notice a possibility of generalizing $\eqref{101}$ to $L$-functions 
which have Euler products and functional equations, because, as explained above, 
the explicit formula is mainly based on the existence of the Euler product and the functional equation. 
However it is not so clear that how we generalize $\eqref{101}$. 
To obtain a hint of the formulation, we recall the outline of the arguments in ~\cite{Su04}. 

A special case of Theorem 1 in ~\cite{Su04} is stated as  
\begin{equation}\label{108}
\sum_{{L(\rho,\chi)=0}\atop{0<{\rm Re}(\rho)<1}} \int_0^\infty h(xu) u^{\rho} \frac{du}{u} =
\sum_{{\zeta(\rho)=0}\atop{0<{\rm Re}(\rho)<1}} \int_0^\infty h(xu) \phi_\chi(u) u^{\rho} \frac{du}{u}
  + O(1) 
\end{equation}
for sufficiently small $x>0$ and $h \in C_0^\infty({\Bbb R}_+)$ with $\int_0^\infty h(u)du=0$, 
where the function $\phi_\chi$ is given by 
\begin{equation}\label{109}
\phi_\chi(u) =\frac{1}{\tau(\bar{\chi})} \sum_{a=1}^q \overline{\chi}(a)e^{2\pi i au/q},
\end{equation}
and $\tau(\chi)=\sum_{a=1}^q \chi(a)e^{2\pi i a/q}$ is the ordinary Gauss sum.  
This is a smooth version of the original result $\eqref{101}$. 
From the well-known equation 
\begin{equation}\label{110}
\chi(n) =\frac{1}{\tau(\bar{\chi})} \sum_{a=1}^q \overline{\chi}(a)e^{2\pi i an/q},
\end{equation}
we find that $\phi_\chi(n)=\chi(n)$ for $n \in {\Bbb Z}$. 
That is, $\phi_\chi$ is an interpolation function of the Dirichlet coefficients $\chi(n)$ of $L(s,\chi)$. 
The existence of such interpolation function plays a key roll in ~\cite{Su04}. 
Now we explain it roughly. 
To obtain $\eqref{108}$, we consider the sum $S(x)=\sum_{n=1}^\infty \Lambda(n)\chi(n)h(xn)$ 
for $h \in C_0^\infty({\Bbb R}_+)$ with $\int_0^\infty h(u)du=0$  
and calculate the sum $S(x)$ in two ways. 
By applying Weil's explicit formula for $L(s,\chi)$ and $u \mapsto h(xu)$,  
we find that $S(x)$ is asymptotically equal to the left hand side of $\eqref{108}$.   
Because $\phi_{\chi}(n)=\chi(n)$, 
we can replace $\chi(n)$ by $\phi_\chi(n)$. 
Denote by $\tilde{S}(x)$ the replaced sum. 
Then, by applying Weil's explicit formula again for $\zeta(s)$ and $u \mapsto h(xu)\phi_{\chi}(u)$,   
we find that $\tilde{S}(x)$ is asymptotically equal to the right hand side of $\eqref{108}$. 
Since $S(x)=\widetilde{S}(x)$, we obtain $\eqref{108}$. 
These arguments suggest that the existence of 
suitable interpolation function of the Dirichlet coefficients 
is very useful for our purpose.  
Standing on the above consideration, 
we adopt 
Euler products, 
functional equations and 
interpolation functions of Dirichlet coefficients 
as the axis of our formulation. 

Here we comment on Euler products. 
In ~\cite{Ku86}, 
Kurokawa showed that the properties of Euler products 
are deeply related to the possibility of analytic continuation. 
His results assert that for a wide class of Euler products 
the unitary property of Euler products is equivalent to 
the possibility of analytic continuation to the whole plane. 
The Euler product also plays a very important role to establish our results. 
Hence our results give a new reason for 
the importance of Euler products in the theory of zeta-functions. 
In our results, the Euler product works as a device which connects the zeros of different $L$-functions with each other. 

This paper is organized as follows. 
In $\S 2$ we prepare several notations and our settings. 
In $\S 3$ we state our theorems. 
In $\S 4$ we give several examples of our results. 
In $\S 5$ we explain about Weil's explicit formula. 
It is a main tool for the proof of our results. 
In $\S 6$ we prove Theorem 1 and Theorem 2. 
In $\S 7$ we prove Theorem 3 and Theorem 4. 
In $\S 8$ we deal with some related topics. 

%
%

\section{ Preliminary. }
In this section we explain our setting and prepare the notations.  
Let $ C^\infty $ be the space of all smooth function on ${\Bbb R}_+ $, and 
let $ C_0^\infty $ be the space of all smooth compactly supported function on ${\Bbb R}_+ $. 
Let $C_D$ be the space of all smooth slowly increasing functions on ${\Bbb R}_+$, that is, 
\begin{equation}\label{201}
C_D:=\{ \phi \in C^\infty  \, | \, |\phi(u)| \leq C u^k \text{ for some $k \geq 0$ and $C>0$ } \}.
\end{equation}
For $\phi \in C_D$ 
we define the Dirichlet series $L_\phi(s)$ as
\begin{equation}\label{202}
L_\phi(s) :=\sum_{n=1}^\infty \phi(n)n^{-s},
\end{equation}
for sufficiently large ${\rm Re}(s)$. 
In the rest of this section, 
we introduce the subclass $C_L$ of $C_D$ that 
any elements in $C_L$ can be regarded  
as an interpolation function of Dirichlet coefficients of some zeta-function. 
To define $C_L$ we recall the concept of the {\it Selberg class}.
The Selberg class $\cal S$ is the class of all function $L(s)$ on ${\Bbb C}$ 
which satisfies the following five axioms;
\begin{enumerate}
\item[(S1)] $L(s)$ is expressed as an absolutely convergent Dirichlet series 
$L(s)=\sum_{n=1}^\infty a(n)n^{-s}$ in the half plane ${\rm Re}(s)>1$.
\item[(S2)] There exists a non-negative integer $m$ such that $(s-1)^m L(s)$ is an entire function with finite order. 
$($We denote by $m_L$ the smallest non-negative integer $m$ which satisfies this condition.$)$
\item[(S3)] 
There exist some $Q>0$, $r \geq 1$, $\lambda_j>0$, ${\rm Re}(\mu_j) \geq 0$ $(1\leq j \leq r)$ and $|\omega|=1$, 
such that {\it the complete $L$-function} 
$$ L^\ast(s):=Q^s \prod_{j=1}^r \Gamma(\lambda_j s + \mu_j) L(s) $$
satisfies the functional equation $L^\ast(s)=\omega \overline{L^\ast(1-\bar{s})}$.
$($ The factor $\gamma(s):=Q^s \prod_{j=1}^r \Gamma(\lambda_j s + \mu_j) $ is called the {\it $\Gamma$-factor} of $L(s)$.$)$  
\item[(S4)] For any positive $\varepsilon$, the coefficients $a(n)$ of $L(s)$ is estimated as $a(n) \ll n^{\varepsilon}$.
\item[(S5)] The logarithm of $L(s)$ is also expressed as  
$$\log L(s)=\sum_{n=1}^\infty b(n)n^{-s},$$ 
where $b(n)$'s are zero unless $n=p^m$ $(m \geq 1)$.  
Further, the estimate $b(n) \ll n^\theta$ holds for some $\theta < \frac 12$.
\end{enumerate}
From (S$5$) the logarithmic derivative of $L(s)$ also has the Dirichlet series expression
\begin{equation}\label{203}
-\frac{L^\prime}{L}(s)
=\sum_{n=1}^\infty \Lambda_L(n)n^{-s}, \quad 
\end{equation}
where $\Lambda_L(n)=b(n)\log n$ is an analogue of 
the von Mangoldt function $\Lambda(n)$ defined by 
\begin{align}\label{204}
\Lambda(n)=
\begin{cases}
\log p & \text{if $n=p^m$ with $m \geq 1$,} \\
0      & \text{otherwise.}
\end{cases}
\end{align}
Futher, by the result of Conrey and Ghosh~\cite{CoGh93}, 
the Dirichlet coefficients $a(n)$ of $L(s) \in {\cal S}$
 are multiplicative. Moreover, the Euler product 
\begin{equation}\label{205}
L(s)=\prod_p L_p(s), \quad \text{where} \quad L_p(s)=\sum_{m=0}^\infty a(p^m) p^{-ms},
\end{equation}
is absolutely convergent for ${\rm Re}(s) >1$ and 
$L_p(s)$ is absolutely convergent for ${\rm Re}(s)>0$ for every $p$. 
The factors $L_p(s)$ are called {\it Euler factors} of $L(s)$. 
We remark that $b(p)=a(p)$ and 
$\log L_p(s)=\sum_{m=1}^\infty b(p^m)p^{-ms}$ for every prime $p$. 
Hence 
\begin{equation}\label{206}
L_p(s) \not= 0 \quad \text{for ${\rm Re}(s)>\theta$ and every $p$,}
\end{equation}
since $b(n) \ll n^\theta$. 
In additon, the $\Gamma$-factor has no zero and no pole 
in the half-plane ${\rm Re}(s) >0$.  
Therefore $L^\ast(s)$ has no zero 
outside of the vertical strip $0 \leq {\rm Re}(s) \leq 1$.  

We say that $L(s) \in {\cal S}$ has a {\it rank $N$ Euler product}, if $L_p(s)$ is expressed as 
\begin{equation}\label{207}
L_p(s) = P_p(p^{-s})^{-1},
\end{equation} 
where
\begin{equation}\label{208}
P_p(X)=1-c_1(p)X-c_2(p)X^2 - \cdots -c_N(p) X^N
\end{equation} 
and $c_N(p) \not=0$ except for finitely many $p$. 
When $L(s)$ has a rank $N$ Euler product, 
we denote by $S_L$ the set of all primes for which $c_N(p)=0$. 

Now we define the subclass $C_L(N)$ and $C_L$ of $C_D$ by    
\begin{align}
C_L(N)  &:= \left\{
 \phi \in C_D \,
\left|
\aligned 
\, & |\phi(u)| \leq C u^2 \quad \text{for some $C>0$,} \label{209} \\
\, & \text{$L_\phi(s)$ belongs to $\cal S$,} \\
\, & \text{$L_\phi(s)$ has a rank $N$ Euler product.}
\endaligned
\right.
\right\}, \\
C_L   &:= \cup_{N \geq 1} C_L(N)\label{210}. 
\end{align}
From the definition, $\phi \in C_L$ satisfies 
the Ramanujan-Deligne estimate $|\phi(n)| \ll_\varepsilon n^{\varepsilon}$ 
for any positive integer $n$, 
even if $\phi(u)$ has rather high order as a function on ${\Bbb R}_+$. 
\smallskip

\noindent
{\bf Remark.} 
The condition $|\phi(u)| \leq C u^2$ is technical one to obtain a simple statement. 
In the case we omit this condition, we can obtain similar results 
although they are of a more complicated form. 
See the proof of Lemma 6 in Section 6.

%
%
\section{Main Results}
\subsection{ Relations with $\zeta(s)$.}
\hfill

\begin{thm}
Let $\phi \in C_L(1)$. 
Then we have
\begin{align} \label{301} \aligned
m_\phi \cdot &  \int_0^\infty h(xu) du  
- \sum_{L_\phi^\ast(\rho)=0} \int_0^\infty h(xu) u^\rho \frac{du}{u} \\
=& 
\, \int_0^\infty h(xu) \phi(u) du  
 - \sum_{\zeta^\ast(\rho)=0} \int_0^\infty h(xu) \phi(u) u^\rho \frac{du}{u} 
+O(1)
\endaligned \end{align}
as $x \to +0$ for any fixed $\varepsilon >0$, 
where $m_\phi=m_{L_\phi}$ is the integer defined in $(${\rm S}$2)$. 
The sum on the left hand side runs over all zeros of $L_\phi^\ast(s)$ counting with multiplicity.
And the sum on the right hand side runs over all zeros of $\zeta^\ast(s)$ counting with multiplicity. 
\end{thm}

\noindent
{\bf Remark.}
Theorem 1 is a generalization of $\eqref{101}$. Compare with~\cite[Theorem 1]{Su04}.\\


\begin{thm}
Let $\phi \in C_L(N)$ with $N \geq 2$. 
Then we have
\begin{align} \label{302} \aligned
m_\phi \cdot &  \int_0^\infty h(xu) du  
- \sum_{L_\phi^\ast(\rho)=0} \int_0^\infty h(xu) u^\rho \frac{du}{u} \\
=& 
\, \int_0^\infty h(xu) \phi(u) du  
 - \sum_{\zeta^\ast(\rho)=0} \int_0^\infty h(xu) \phi(u) u^\rho \frac{du}{u} 
+ J_{\phi}(x;h) + O(x^{-\frac 13 -\varepsilon})
\endaligned \end{align}
as $x \to +0$ for any fixed $\varepsilon >0$, 
where $m_\phi=m_{L_\phi}$ is the integer defined in $(${\rm S}$2)$. 
The sum on the left hand side runs over all zeros of $L_\phi^\ast(s)$ counting with multiplicity, 
and the sum on the right hand side runs over all zeros of $\zeta^\ast(s)$ counting with multiplicity. 
The function $J_{\phi}(x;h)$ is estimated as 
\begin{equation}\label{303}
J_{\phi}(x;h) \ll_{\varepsilon, \phi, h} x^{-\frac{1}{2}-\varepsilon} 
\end{equation}
as $x \to +0$ for any fixed $\varepsilon >0$. 
Furthermore, if the estimate 
\begin{equation}\label{304}
\sum_{p \leq T} c_{\phi,2}(p) \log p = A_\phi T^{\mu}  +O(T^{\nu +\varepsilon}).
\end{equation}
holds for the numbers $c_{\phi,2}(p)$ defined in $(2.8)$ with some constant $A_\phi$ and ${\rm Re}(\mu) > \nu \geq 0$, 
then 
\begin{equation}\label{305}
J_{\phi}(x;h) = \mu A_\phi \int_0^\infty h(xu^2) u^{\mu} \frac{du}{u}  + O(x^{-\frac {\nu}{2} -\varepsilon}).
\end{equation}
\end{thm}

\noindent
{\bf Remark.}
When $A_{\phi} \not= 0$, $J_{\phi}(x;h) \sim C_{\phi,h}x^{-\frac{\mu}{2}}$. 
Hence ${\rm Re}(\mu)$ must be smaller than $1$ by $\eqref{303}$ in this case. 
When $A_{\phi}=0$, $J_{\phi}(x;h)=O(x^{-\frac{\nu}{2} -\varepsilon})$.  
Hence $\nu$ must be smaller than $1$ by $\eqref{303}$ in this case. 


\begin{cor} 
Let $\phi, \,\psi \in C_L(N)$. Suppose that $\phi(n)=\psi(n)$ for any positive integer $n$. 
In the case $N \geq 2$, 
we additionally suppose that $\phi$ $($and $\psi)$ 
satisfies the condition $(3.4)$ for the constants $A,\mu,\nu$. 
Then we have 
\begin{align}\label{306}
\aligned
\int_0^\infty & h(xu)\phi(u)du 
- \sum_{\zeta^\ast(\rho)=0} \int_0^\infty h(xu) \phi(u) u^\rho \frac{du}{u} \\
=& 
\int_0^\infty h(xu) \psi(u)du 
- \sum_{\zeta^\ast(\rho)=0} \int_0^\infty h(xu) \psi(u) u^\rho \frac{du}{u} 
+O(x^{-\frac {\nu}{2} -\varepsilon}) +O(x^{-\frac 13 -\varepsilon})
\endaligned
\end{align}
as $x \to 0$ for any fixed $\varepsilon >0$. 
Sums on the both sides run over all zeros of $\zeta^\ast(s)$ counting with multiplicity. 
\end{cor}

\subsection{ Rankin-Selberg type relations .} 
Let $\phi \in C_L(M)$, $\psi \in C_L(N)$ with the Euler products 
\begin{align*}
L_\phi(s)& =\prod_p P_{\phi,p}(p^{-s})^{-1} 
      =\prod_{p} \prod_{i=1}^{{\rm deg}\,P_{\phi,p}} (1 -\alpha_{\phi}(p,i) p^{-s})^{-1}, \\
L_\psi(s)& =\prod_p P_{\psi,p}(p^{-s})^{-1} 
      =\prod_{p} \prod_{j=1}^{{\rm deg}\,P_{\psi,p}} (1-\alpha_{\psi}(p,j) p^{-s})^{-1}.
\end{align*}
We use the notaion $S=S_{\phi \otimes \psi}=S_{\phi} \cup S_{\psi}$ 
and denote by $c_{\phi,l}(p)$ the $l$-th coefficient of polynomial $P_{\phi,p}(X)$ (cf. $\eqref{208}$).  
Define the function $\widetilde{L}_{\psi \otimes \phi}(s)$ by  
\begin{equation}\label{307}
\widetilde{L}_{\phi \otimes \psi}(s)
=\prod_{p \not\in S} \prod_{i=1}^M \prod_{j=1}^N (1-\alpha_{\phi}(p,i) \alpha_{\psi}(p,j)p^{-s})^{-1}
\end{equation}
for sufficiently large ${\rm Re}(s)$. Note that 
\begin{equation}\label{308}
\widetilde{L}_{\phi \otimes \psi}(s)
= \prod_{p \not\in S} \prod_{i}(1-\alpha_{\phi}(p,j) \psi(p)p^{-s})^{-1} 
= \prod_{p \not\in S} P_{\phi,p}( \psi(p)p^{-s} )^{-1}, 
\end{equation}
if $\psi \in C_L(1)$. 
We say that the pair $(\phi,\psi) \in C_L(M) \times C_L(N)$ is the {\it Rankin-Selberg pair}, 
if there exist some polynomials $Q_p(X)$ with ${\rm deg}\,Q_p \leq MN$ for any $p \in S$ 
such that the function
\begin{equation}\label{309}
L_{\phi \otimes \psi}(s):=
\widetilde{L}_{\phi \otimes \psi}(s) \times 
\prod_{p \in S} Q_p(p^{-s})^{-1}
\end{equation}
belongs to the Selberg class $S$. 

\begin{thm} 
Let $(\phi,\psi) \in C_L(N) \times C_L(1)$ be a Rankin-Selberg pair. 
Then we have
\begin{equation}\label{310a}
\aligned
m_{\phi \otimes \psi} \cdot  & \int_0^\infty h(xu) du  
- \sum_{L_{\phi \otimes \psi}^\ast(\rho)=0} \int_0^\infty h(xu) u^\rho \frac{du}{u} \\
=& 
\,m_\phi \cdot \int_0^\infty h(xu) \psi(u) du  
- \sum_{L_\phi^\ast(\rho)=0} \int_0^\infty h(xu) \psi(u) u^\rho \frac{du}{u}
+ O(x^{-\frac{1}{3}-\varepsilon}) 
\endaligned
\end{equation}
as $x \to 0$ for any fixed $\varepsilon$. 
Additionally, we suppose that $|\phi(u)\psi(u)|\leq C u^2$ for some $C>0$. 
Then 
\begin{equation}\label{310b}
\aligned
m_{\phi \otimes \psi} \cdot  & \int_0^\infty h(xu) du  
- \sum_{L_{\phi \otimes \psi}^\ast(\rho)=0} \int_0^\infty h(xu) u^\rho \frac{du}{u} \\
=& 
\,m_\phi \cdot \int_0^\infty h(xu) \psi(u) du  
- \sum_{L_\phi^\ast(\rho)=0} \int_0^\infty h(xu) \psi(u) u^\rho \frac{du}{u}
+ O(x^{-\frac{1}{3}-\varepsilon}) \\
=& 
\, \int_0^\infty h(xu) \phi(u) \psi(u) du  
- \sum_{\zeta^\ast(\rho)=0} \int_0^\infty h(xu) \phi(u) \psi(u) u^\rho \frac{du}{u}
+ J_{\phi,\psi}(x;h)  + O(x^{-\frac{1}{3}-\varepsilon})
\endaligned
\end{equation}
as $x \to 0$ for any fixed $\varepsilon$. 
The function $J_{\phi,\psi}(x;h)$ is estimated as 
\begin{equation}\label{311}
J_{\phi,\psi}(x;h) \ll_{\varepsilon, \phi,\psi, h} x^{-\frac{1}{2}-\varepsilon} 
\end{equation}
as $x \to +0$ for any fixed $\varepsilon >0$. 
Furthermore, if the estimate 
\begin{equation}\label{312}
\sum_{p \leq T} c_{\phi,2}(p) \psi(p)^2 \log p = A_{\phi,\psi} T^{\mu}  +O(T^{\nu +\varepsilon}).
\end{equation}
holds with some constant $A_{\phi,\psi}$ and ${\rm Re}(\mu) > \nu \geq 0$, 
\begin{equation}\label{313}
J_{\phi,\psi}(x;h) = \mu A_{\phi,\psi} \int_0^\infty h(xu^2) u^{\mu} \frac{du}{u}  + O(x^{-\frac {\nu}{2} -\varepsilon}).
\end{equation}
\end{thm}


\begin{thm} 
Let $(\phi,\psi) \in C_L(M) \times C_L(N)$ be a Rankin-Selberg pair with $M,N \geq 2$. 
Then 
\begin{equation}\label{314a}
\aligned
m_{\phi \otimes \psi} & \int_0^\infty h(xu) du  
- \sum_{L_{\phi \otimes \psi}^\ast(\rho)=0} \int_0^\infty h(xu) u^\rho \frac{du}{u} \\
=& 
\,m_\phi \int_0^\infty h(xu) \mu(u) du  
  - \sum_{L_\phi^\ast(\rho)=0} \int_0^\infty h(xu) \psi(u) u^\rho \frac{du}{u} \\
& +J_{\phi,\psi}^{(1)}(x;h) +J_{\phi,\psi}^{(2)}(x;h)  + O(x^{-\frac{1}{3}-\varepsilon })  \\ 
=& 
\, \int_0^\infty h(xu) \phi(u) \psi(u) du  
- \sum_{\zeta^\ast(\rho)=0} \int_0^\infty h(xu) \phi(u) \psi(u) u^\rho \frac{du}{u} \\
& +J_{\phi,\psi}^{(3)}(x;h) + J_{\phi,\psi}^{(4)}(x;h) + J_{\phi,\psi}^{(5)}(x;h) + O(x^{-\frac{1}{3}-\varepsilon })
\endaligned
\end{equation}
as $x \to +0$ for any fixed $\varepsilon >0$. 
Additionally, we suppose that $|\phi(u)\psi(u)|\leq C u^2$ for some $C>0$. 
Then 
\begin{equation}\label{314b}
\aligned
m_{\phi \otimes \psi} & \int_0^\infty h(xu) du  
- \sum_{L_{\phi \otimes \psi}^\ast(\rho)=0} \int_0^\infty h(xu) u^\rho \frac{du}{u} \\
=& 
\,m_\phi \int_0^\infty h(xu) \mu(u) du  
  - \sum_{L_\phi^\ast(\rho)=0} \int_0^\infty h(xu) \psi(u) u^\rho \frac{du}{u} \\
& +J_{\phi,\psi}^{(1)}(x;h) +J_{\phi,\psi}^{(2)}(x;h)  + O(x^{-\frac{1}{3}-\varepsilon })  \\ 
=& 
\, \int_0^\infty h(xu) \phi(u) \psi(u) du  
- \sum_{\zeta^\ast(\rho)=0} \int_0^\infty h(xu) \phi(u) \psi(u) u^\rho \frac{du}{u} \\
& +J_{\phi,\psi}^{(3)}(x;h) + J_{\phi,\psi}^{(4)}(x;h) + J_{\phi,\psi}^{(5)}(x;h) + O(x^{-\frac{1}{3}-\varepsilon })
\endaligned
\end{equation}
as $x \to +0$ for any fixed $\varepsilon >0$. 
The functions $ J_{\phi,\psi}^{(k)}(x;h) $ are estimated as 
\begin{equation}\label{315}
J_{\phi,\psi}^{(k)}(x;h) \ll x^{-\frac{1}{2}-\varepsilon } \quad (1 \leq k \leq 5)
\end{equation} 
as $x \to +0$ for any fixed $\varepsilon >0$. 
Moreover the asymptotic formulas  
\begin{align}
J_{\phi,\psi}^{(1)}(x;h) = J_{\phi,\psi}^{(3)}(x;h) 
&= \mu_1 A_1 \int_0^\infty h(xu^2) u^{\mu_1} \frac{du}{u}
 +O(x^{-\frac{\nu_1}{2}-\varepsilon}),  \label{316} \\
J_{\phi,\psi}^{(4)}(x;h) 
&=  \mu_2 A_2 \int_0^\infty h(xu^2) u^{\mu_2} \frac{du}{u}
 +O(x^{-\frac{\nu_2}{2} -\varepsilon}),  \label{317} \\
\frac{1}{2}J_{\phi,\psi}^{(2)}(x;h) = \frac{1}{3}J_{\phi,\psi}^{(5)}(x;h)
&= \mu_3 A_3 \int_0^\infty h(xu^2) u^{\mu_3} \frac{du}{u}
 +O(x^{-\frac{\nu_3}{2} -\varepsilon}) \label{318}
\end{align}
hold, if the corresponding estimate
\begin{align}
\sum_{p \leq T}\Lambda_{\phi \otimes \phi}(p) c_{\psi,2}(p) = A_1 T^{\mu_1}+O(T^{ \nu_1+\varepsilon }), \label{thm4:1} \\
\sum_{p \leq T} \Lambda_{\psi \otimes \psi}(p) c_{\phi,2}(p) = A_2 T^{\mu_2}+O(T^{\nu_2+\varepsilon}),\label{thm4:2} \\
\sum_{p \leq T} c_{\phi,2}(p) c_{\psi,2}(p) \log p = A_3 T^{\mu_3}+O(T^{ \nu_3+\varepsilon }),\label{thm4:3} 
\end{align}
hold respectively.
\end{thm}

%
%
\section{Examples}

In this section we give simple examples of Theorem 1 $\sim$ Theorem 4. 

%
%

\subsection{ Dirichlet $L$-functions. } 
Let $\chi$ mod $q$ be a primitive Dirichlet character.  
Then the function $\phi_{\chi}$ defined in $\eqref{109}$ belongs to $C_L(1)$, 
since $\phi_\chi$ is bounded on ${\Bbb R}_+$ and $L_{\phi_{\chi}}(s)=L(s,\chi)$. 
Further 
$$
\int_0^\infty h(xu) \phi_{\chi}(u)du
= \frac{-q}{2\pi i \,\tau(\overline{\chi})} \sum_{a=1}^q \frac{\overline{\chi}(a)}{a}
\int_0^\infty h^\prime(v)e^{2 \pi i \frac{av}{qx}}dv =O(1).
$$ 
Hence we re-obtain $\eqref{108}$ from Theorem 1. 

%
%

\subsection{ Automorphic $L$-functions attached to cusp forms in $S_k(N)$. } 
Let ${\frak h}=\{ z \in{\Bbb C};{\rm Im}(z)>0\}$ be the upper half plane 
and let $\Gamma_0(N)$ be the Hecke subgroup of level $N$ of the full modular group.   
Let $S_k(N)$ be the vector space of all holomorhic function $f$ on $\frak h$ such that 
$f((az+b)/(cz+d))=(cz+d)^k f(z)$ 
for any $\bigl( \begin{smallmatrix} a&b\\c&d \end{smallmatrix} \bigr) \in \Gamma_0(N)$,  
and $f(i\infty)=0$. It is well known that any $f \in S_k(N)$ has the Fourier expansion 
\begin{equation}\label{401}
f(z)=\sum_{n=1}^\infty a_{f}(n) e^{ 2 \pi i n z}.
\end{equation}
By using the Fourier coefficients $\{a_{f}(n) \}$, 
the automorpchic $L$-function $L(s,f)$ associated with $f$ is defined as 
\begin{equation}\label{402}
L(s,f)=\sum_{n=1}^\infty a_{f}(n) n^{-s-\frac{k-1}{2}}.
\end{equation}
This series is absolutely convergent on the right-half plane ${\rm Re }(s)>1$ 
because of the estimate $\sum_{n \leq T} |a_{f}(n)|^2 \ll T^{k+1}$ 
obtained by the Rankin-Selberg method or 
the more precise estimate $|a_{f}(n)| \ll_\varepsilon n^{\varepsilon}$ due to Deligne.
The automorphic $L$-function $L(s,f)$ can be extended to an entire function in $s$ and the function
$$L^\ast(s,f)=N^{\frac{2s+k-1}{4}}(2\pi)^{-s-\frac{k-1}{2}} \Gamma(s+\frac{k-1}{2}) L(s,f)$$
satisfies the functional equation 
$$L^\ast(s,f)=\pm (-1)^{k/2} L^\ast(1-s,f)$$ 
where the sign $\pm$ is determined by the action of Fricke involution.
Moreover, if $f\in S_k(N)$ is a normalized Hecke-eigen newform (cf. ~\cite[chap.1.4]{Bu99}), 
$L(s,f)$ has the Euler product 
\begin{equation}\label{403}
L(s,f)=
\prod_{p|N} (1-a_{f}(p)p^{-\frac{k-1}{2}}p^{-s})^{-1}
\prod_{p \not\,|\, N} (1-a_{f}(p)p^{-\frac{k-1}{2}}p^{-s}+p^{-2s})^{-1}.
\end{equation}
Define the function $\phi_f:{\Bbb R}_+  \to {\Bbb C}$ by 
\begin{equation}\label{404}
\phi_f(u)=u^{-\frac{k-1}{2}}\,\int_1^2 f(X+iu^{-1})e^{ -2\pi i u (X+iu^{-1}) }dX. \label{Def_1}
\end{equation}
From the definiton of $\phi_f$, $\phi_f(n)$ coincides with the (shifted) $n$-th Fourier coefficient $a_{f}(n)n^{-\frac{k-1}{2}}$,  
and satisfies the estimate $|\phi_f(u)| \leq C\sqrt{u}$ for some $C>0$, because $Y^{k/2}|f(X+iY)|$ is bounded.  
Also we can easily find that $c_2(p)=-1$ and $S_{\phi_f}=\{p\,;\, p|N \}$ for $\phi_f(u)$. 
Therefore $\phi_f(u)$ belongs to $C_L(2)$ and 
we find that 
\begin{equation}\label{405}
\int_0^\infty h(xu)\phi_f(u)du=O(x^{N-\frac{3}{2}})
\end{equation}
for any fixed positive integer $N$ by using integration by parts suitable times. 
Furthermore, if we assume that $\zeta(s)$ has no zero in ${\rm Re}(s) > \sigma$, 
then  
\begin{equation}\label{406}
\sum_{p\leq T, \,p \not\in S_{\phi} }c_2(p)\log p 
= \sum_{p\leq T, \,p \not\in S_{\phi} } \log p = T +O(T^{\sigma+\varepsilon}).
\end{equation}
Hence we obtain the following result as a consequence of Theorem 2.


\begin{thm} 
Let $f \in S_k(N)$ be a normalized Hecke-eigen cuspform. 
Assume that $\zeta(s)$ has no zeros in ${\rm Re}(s) > \sigma$. 
Then
\begin{equation}\label{407}
\aligned
\sum_{L^\ast(\rho,f)=0} & \int_0^\infty 
h(xu)u^{\rho}\frac{du}{u} \\
= & \sum_{\zeta^\ast(\rho)=0}\int_0^\infty 
h(xu) \phi_f(u) u^\rho \frac{du}{u}
-C(h) \, x^{-\frac 12} 
+O(x^{-\frac{\sigma}{2}-\varepsilon }) +O(x^{-\frac 13 - \varepsilon }) 
\endaligned 
\end{equation}
as $x \to +0$ for any fixed $\epsilon >0$ for any $h \in C_0^\infty$, 
where $C(h)=2^{-1}\,\widehat{h}(1/2)$. 
The sum on the right hand side runs over all zeros of $\zeta^\ast(s)$ counting with multiplicity. 
And the sum on the left hand side runs over all zeros of $L^\ast(s,f)$ counting with multiplicity.
\end{thm}

%
%

\subsection{ Rankin-Selberg $L$-functions. }
Let $f,g \in S_k(1)$ be normalized Hecke-eigen cusp forms with Fourier expansions 
\begin{equation}\label{408}
f(z)=\sum_{n=1}^\infty a_{f}(n)e^{ 2 \pi i n z},\qquad
g(z)=\sum_{n=1}^\infty a_{g}(n)e^{ 2 \pi i n z}.
\end{equation}
Define $\alpha_p$, $\beta_p$, $\gamma_p$ and $\delta_p$ 
by using the Euler product of $L(s,f)$, $L(s,g)$; 
\begin{align*}
L(s,f)&=\prod_p (1 -a_f(p)p^{-\frac{k-1}{2}} p^{-s} +p^{-2s})^{-1} 
      =\prod_p [(1 -\alpha_p p^{-s}) (1-\beta_p p^{-s})]^{-1}, \\
L(s,g)&=\prod_p (1 -a_g(p)p^{-\frac{k-1}{2}} p^{-s} +p^{-2s})^{-1} 
      =\prod_p [(1-\gamma_p p^{-s})(1-\delta_p p^{-s})]^{-1}.
\end{align*}
The Rankin-Selberg $L$-function $L(s,f \otimes g)$ is defined by 
\begin{equation*}
L(s,f \otimes g)=\prod_p 
[(1-\alpha_p \gamma_p p^{-s})(1-\alpha_p \delta_p p^{-s})(1-\beta_p \gamma_p p^{-s})(1-\beta_p \delta_p p^{-s})]^{-1}.
\end{equation*}
Then $L(s,f \otimes g)=\zeta(2s)\sum_{n=1}^\infty a_{f}(n) a_{g}(n) n^{-s-k+1}$.
Moreover the completed $L$-function $L^\ast(s, f \otimes g)=(4\pi)^{-s-k+1}\Gamma(s+k-1)\Gamma(s)L(s, f \otimes g)$ 
satisfies the functional equation 
$L^\ast(s, f \otimes g)=L^\ast(1-s, f \otimes g)$ 
(cf. ~\cite[chap.1.6]{Bu99}). 
Further it is known that $s=1$ is a simple pole if $f=g$ and is a regular point otherwise. 
Let 
\begin{equation}\label{409}
\phi_{f \otimes g}(u):= \phi_f(u)\phi_g(u)
\end{equation}
where $\phi_f$, $\phi_g$ are defined in $\eqref{404}$. 
Then $ |\phi_{f \otimes g}(u)| \leq C u $ for some $C>0$. 
If $\zeta(s)$, $L(s,f \otimes f)$ and $L(s,g \otimes g)$ have no zeros in ${\rm Re}(s)>\sigma$, 
we obtain
\begin{equation}\label{410}
\aligned 
\sum_{p \leq T}\Lambda_{\phi_f \otimes \phi_f}(p) c_{\phi_g,2}(p) 
= -\sum_{p \leq T}\Lambda_{f \otimes f}(p) = - T +O(T^{ \sigma+\varepsilon }),  \\
\sum_{p \leq T}\Lambda_{\phi_g \otimes \phi_g}(p) c_{\phi_g,2}(p) 
= -\sum_{p \leq T} \Lambda_{g \otimes g}(p) = - T +O(T^{\sigma+\varepsilon}), \\
\sum_{p \leq T} c_{\phi_f,2}(p) c_{\phi_g,2}(p) \log p 
=\sum_{p \leq T} \log p = T +O(T^{ \sigma+\varepsilon }). 
\endaligned \end{equation}
Together with the above things, we obtain the following theorem as a consequence of Theorem 4. 


\begin{thm} 
Let $f,g \in S_k(1)$ be normalized Hecke-eigen cusp forms. 
Assume that $\zeta(s)$, $L(s,f \otimes f)$ and $L(s,g \otimes g)$
have no zeros in ${\rm Re}(s) >\sigma$. 
Then for any $h \in C_0^\infty$, the following formula holds:
\begin{equation}\label{411}
\aligned
\delta_{f=g} & \int_{0}^\infty h(xu)du 
-\sum_{L^\ast(\rho,f \otimes g)=0}\int_0^\infty 
h(xu) u^{\rho} \frac{du}{u} \\
\,
&= -\sum_{L^\ast(\rho,f)=0} \int_0^\infty 
h(xu) \phi_g(u) u^\rho \frac{du}{u}
- C(h) \, x^{-\frac 12}
+O(x^{-\frac {\sigma}{2} - \varepsilon }) 
+O(x^{-\frac 13 - \varepsilon }) \\
\,
&= -\sum_{\zeta^\ast(\rho)=0}\int_0^\infty 
h(xu) \phi_f(u) \phi_g(u) u^\rho \frac{du}{u}
- C(h) \, x^{-\frac 12}
+O(x^{-\frac {\sigma}{2} - \varepsilon }) 
+O(x^{-\frac 13 - \varepsilon } ) 
\endaligned 
\end{equation}
as $x \to +0$ for any positive $\varepsilon$, 
where $ C(h) =2^{-1} \, \widehat{h}(1/2)$, $\delta_{f=g}=1$ if $f=g$ and is zero otherwise.
\end{thm}

%
%
\section{Weil's Explicit Formula}
In this section we state a version of Weil's explicit formula. 
It is one of the main tools for our proof of the results in this paper. 
Define the involution $h \mapsto h^\ast$ on $C_0^\infty$ by 
\begin{equation}\label{501}
h^\ast(u)= \frac 1u f(\frac 1u)
\end{equation}
and the Mellin transform of $h$ by 
\begin{equation}\label{502}
\widehat{h}(s)=\int_0^\infty h(u)u^s \frac{du}{u}.
\end{equation}
Because $h$ has a compact support, the above integral is absolutely convergent for any $s \in {\Bbb C}$. 
Further the Mellin inversion formula 
\begin{equation}\label{503}
h(u)=\frac{1}{2\pi i}\int_{(\sigma)}\widehat{h}(s)u^{-s}ds
\end{equation}
is valid, where the path of integration is the vertical line ${\rm Re}(s)=\sigma$. 


\begin{prop} {\bf [ Weil's Explicit Formula ]} 
Let $L(s) \in {\cal S}$. 
Then, for any $h \in C_0^\infty$, 
\begin{equation}\label{504}
\aligned
m_L & \, \widehat{h}(0) -\sum_{L^\ast(\rho)=0}\widehat{h}(\rho) + m_L \, \widehat{h}(1) \\
=&\sum_{n=1}^\infty \{ \Lambda_L(n) h(n)+\overline{\Lambda_L(n)} h^\ast (n) \}
  +(2 \log Q+d\,C_E)\,h(1)+\sum_{j=1}^r W_{\lambda_j.\mu_j}(h),
\endaligned 
\end{equation}
where $d=2\sum_{j=1}^r\lambda_j$ and $C_E$ is the Euler constant. 
The functional $W_{\lambda,\mu}$ is given by 
\begin{equation}\label{505}
W_{\lambda,\mu}(h) = 
\int_1^\infty 
\left[h_{\lambda,\mu}(u)+h_{\lambda,\mu}^\ast(u)-2h(1)u^{({\rm Re}(\mu)-1)/\lambda}\right]
\frac{u^{(1-{\rm Re}(\mu))/\lambda}}{u^{1/\lambda}-1}\frac{du}{u}, 
\end{equation}
\begin{equation}\label{506}
h_{\lambda,\mu}(u)  =h(u)u^{-i\,{\rm Im}(\mu)/\lambda}. 
\end{equation} 
The sum $\sum_{L^\ast(\rho)=0}$ runs over all zeros of $L^\ast(s)$ counting with multiplicity. 
Sums and integrals contained in the both sides of $\eqref{504}$ are absolutely convergent, 
because the Mellin transform $\widehat{h}$ decays very fast by the assumption on $h$.
\end{prop}

Proposition 1 is proved by a way similar to the proof of Weil's explicit formula in ~\cite{La70}. 
There is no essential difference or difficulty in our case 
because of conditions $({\rm S}1) \sim ({\rm S}5)$ for $L(s)$. 
Hence we omit the proof of Proposition 1.

%
%
\section{Proof of Theorem 1 and Theorem 2}

\subsection{ Lemmas.} 
In this part we prepare several lemmas which are necessary for our proof of Theorem 1 and Theorem2. 
For $L(s) \in {\cal S}$, $h \in C_0^\infty$ and $x>0$, 
we define the sum $S_L(x)$ by 
\begin{equation}\label{601}
S_L(x) :=S_L(x;h) :=\sum_{n=1}^\infty \Lambda_L(n)h(xn).
\end{equation}


\begin{lem}
Let $L(s)=\sum_{n=1}^\infty c(n) n^{-s} \in {\cal S}$. 
Then, for any $\varepsilon >0$, 
\begin{equation}\label{602}
\sum_{p} \sum_{{l \leq m}\atop{p^m \leq T} } c(p^m) \log p
\ll_\varepsilon  T^{\frac 1l +\varepsilon } ,
\end{equation}
and
\begin{equation}\label{603}
\sum_{p} \sum_{{l \leq m}\atop{p^m \leq T} } \Lambda_L(p^m)
\ll_\varepsilon  T^{ \theta +\frac 1l +\varepsilon } ,
\end{equation} 
where $\theta$ is the constant in axiom $({\rm S}5)$ of the Selberg class. 
\end{lem}


\begin{pf} 
We have 
\begin{equation}\label{604}
\aligned
\sum_{p} &\sum_{{l \leq m}\atop{p^m \leq T} } c(p^m) \log p
  \ll \sum_{p^l \leq T} \sum_{ l \leq m \leq \frac{\log T}{\log p} } p^{m\varepsilon} \log p \\
& \leq T^\varepsilon \frac 1l \log T \sum_{p^l \leq T} \sum_{ 1 \leq m \leq \frac{\log T}{\log p} } 1 
  \leq \frac {1}{\log2} T^{\varepsilon +\frac 1l} (\log T)^2 
  \ll_{\varepsilon^\prime} T^{\frac 1l +\varepsilon^\prime}.
\endaligned
\end{equation}
This is assertion $\eqref{602}$.
Recall that $\Lambda_L(n) = b(n)\log n$ and $|b(n)| \ll n^\theta$. 
Then we have
\begin{equation}\label{605}
\aligned
\sum_{p} & \sum_{{l \leq m}\atop{p^m \leq T} } \Lambda_L(p^m) 
  = \sum_{p^l \leq T} \sum_{ l \leq m \leq \frac{\log T}{\log p} } b(p^m)\log p^m \\
& \ll \sum_{p^l \leq T} \sum_{ l \leq m \leq \frac{\log T}{\log p} } p^{m\theta} \log p^m 
  \leq T^\theta \log T \sum_{p^l \leq T} \sum_{ 1 \leq m \leq \frac{\log T}{\log p} } 1 \\
& \leq T^{\theta + \frac 1l} \log T \frac {\log T}{\log 2} \sum_{p^l \leq T} 
  \ll_\varepsilon T^{\theta +\frac 1l +\varepsilon}. 
\endaligned
\end{equation}
This is our assertion $\eqref{603}$. 
\end{pf}


\begin{lem} 
Let $L(s) \in {\cal S}$. Then
\begin{equation}\label{606}
S_L(x)
=m_L \cdot \int_0^\infty h(xu) du 
-\sum_{L^\ast(\rho)=0} \int_0^\infty h(xu) u^\rho \frac{du}{u}
+O(1).
\end{equation}
\end{lem}

\begin{pf}
Applying Proposition 1 to $L(s)$ and $u \mapsto h(xu)$ we have 
\begin{equation}\label{607} 
\aligned
S_L(x)=\,\,
& m_L \int_0^\infty h(xu) du
 -\sum_{L^\ast(\rho)=0} \int h(xu)u^\rho \frac{du}{u} 
 + m_L \int_0^\infty h(xu) \frac{du}{u} \\
& -\sum_{n=1}^\infty  \overline{\Lambda_L(n)} n^{-1} h(xn^{-1}) \\
& -(2 \log Q + d\,C_E)\,h(x)
 -\sum_{j=1}^r W_{\lambda_j,\mu_j}(u \mapsto h(xu)).
\endaligned 
\end{equation}

The third term on the right hand side is bounded, because it is equal to $\widehat{h}(0)$. 
The fourth term and the fifth term on the right hand side is zero 
for sufficiently small $x>0$ because the support of $h$ is compact. 
Furthermore the sixth term on the right hand side is absorbed into the error term. 
In fact 
$$
W_{\lambda, \mu}(u \mapsto h(xu)) 
= \int_1^\infty h(xu) 
\frac{u^{\frac{1-\mu}{\lambda}}}{u^{ \frac{1}{\lambda} } -1}\frac{du}{u} 
= x^{\frac{\mu}{\lambda}} 
\int_x^\infty h(v) 
\frac{v^{ \frac{1-\mu}{\lambda} }}{ v^{ \frac{1}{\lambda} } -x^{\frac{1}{\lambda} }}\frac{dv}{v}
$$
for sufficiently small $x>0$ and the right hand side is bounded as $x \to +0$ 
since ${\rm Re}(\mu) \geq 0$ and $\lambda >0$.
\end{pf}


\begin{lem} Let $L(s) \in {\cal S}$. Then
\begin{equation}\label{608}
S_L(x)
=\sum_{p} \sum_{m=1}^{l-1} \Lambda_L(p^m) h( x p^m )
+O(x^{ -\theta -\frac 1l -\varepsilon})
\end{equation}
for any $\varepsilon >0$.
\end{lem}

\begin{pf}
It suffices to show that 
\begin{equation}\label{609}
\sum_{p} \sum_{m=l}^{\infty} \Lambda_L(p^m) h( x p^m ) = O(x^{ -\theta -\frac 1l -\varepsilon}).
\end{equation}
By using partial summation (cf.~\cite[page 2]{SP94}) we obtain 
$$
\sum_{p} \sum_{m=l}^{\infty} \Lambda_L(p^m) h( x p^m ) 
= - \int_l^\infty xu h^\prime (xu) \left( \sum_{p}\sum_{{l \leq m}\atop{p^m \leq u}} \Lambda_L(p^m) \right) \frac{du}{u}.
$$
From Lemma 1 the right hand side is estimated as 
$$
O\left( \int_{l}^\infty
xu |h^\prime(xu)|u^{ \theta +\frac{1}{l} +\varepsilon }\frac{du}{u}
\right)
=O( x^{ -\theta -\frac{1}{l} -\varepsilon } ).
$$
This implies $\eqref{609}$.
\end{pf}

Here we describe the relation between the Dirichlet coefficients 
of $L_\phi(s)$ and those of its logarithmic derivative $(L_\phi^\prime/L_\phi)(s)$. 
Define the numbers $r_m(p)$ by 
\begin{equation}\label{610}
X \frac{d}{dX} \log \, P_p(X)^{-1} = \sum_{m=1}^\infty r_m(p)X^m.
\end{equation}
By simple series calculations, we find that 
\begin{equation}\label{611}
\phi(p^m)= 
\begin{cases}
c_1(p)\phi(p^{m-1})+c_2(p)\phi(p^{m-2})+\cdots+c_m(p) & \text{if $m \leq n_p$} \\
c_1(p)\phi(p^{m-1})+c_2(p)\phi(p^{m-2})+\cdots+c_{n_p}(p)\phi(p^{m-n_p}) & \text{if $m > n_p$}
\end{cases}
\end{equation}
and 
\begin{equation}\label{612}
r_m(p)= \left\{
\aligned
\phi( & p^m)  +c_2(p)\phi(p^{m-2}) +\cdots \\
& +(j-1)c_j(p)\phi(p^{m-j}) +\cdots +(m-1)c_m(p) \qquad \quad \quad \,\,\,\,\,\,\, \text{if $m \leq n_p$} \\
\phi( & p^m)  +c_2(p)\phi(p^{m-2}) +\cdots \\
& +(j-1)c_j(p)\phi(p^{m-j}) +\cdots +(n_p-1)c_{n_p}(p)\phi(p^{m-n_p}) \quad  \text{if $m > n_p$} 
\endaligned 
\right.
\end{equation}
where $n_p={\rm deg}\,P_p(X)$ and 
the numbers $c_j(p)$ are the coefficients of the polynomial $P_p(X)$ 
defined in $\eqref{207}$, $\eqref{208}$. 
Note that $\phi(p^0)=\phi(1)=1$ 
which can be seen from the form of the Euler product attached to $L_\phi(s)$. 
Then the Dirichlet coefficient $\Lambda_\phi(n)$ of $(L_\phi^\prime/L_\phi)(s)$ is given by
\begin{equation}\label{613}
\Lambda_{\phi}(n) =
\begin{cases} 
r_m(p) \log p & \text{ if $n=p^m$ with $m \geq 1$}, \\
0             & \text{otherwise.}
\end{cases}
\end{equation}
Additionally, it is useful for us to note the relation 
\begin{equation}\label{614}
c_l(p) = (-1)^{l+1} \sum_{ {(i_1,\cdots,i_l)}\atop{1\leq i_1 \leq \cdots \leq i_l \leq M} } 
\alpha_{\phi}(p,i_1) \cdots \alpha_{\phi}(p,i_l) 
\end{equation}
where $\alpha_{\phi}(p,i) $ are the roots of 
the polynomial $P_{\phi,p}(X)$ associated with $\phi \in C_L(M)$.


\begin{lem}
Let $\phi \in C_L(N)$.
Then 
\begin{equation}\label{615}
\Lambda_{\phi}(n) \ll_\varepsilon n^{\varepsilon}. 
\end{equation}
\end{lem}


\begin{pf} 
From $\eqref{611}$ we have 
\begin{equation}\label{616}
c_j(p) = \sum_{l=1}^j \sum_{{1\leq i_1,\cdots,i_l \leq j}\atop{i_1+\cdots+i_l=j}} 
(-1)^{l+1}\phi(p^{i_1}) \cdots \phi(p^{i_l}).
\end{equation}
Because $|\phi(n)| \leq C_\varepsilon n^\varepsilon$ for any positive integer $n$, 
$$
|c_j(p)| \leq  
p^{j \varepsilon}
\sum_{l=1}^j \sum_{{1\leq i_1,\cdots,i_l \leq j}\atop{i_1+\cdots+i_l=j}} 
C_\varepsilon^l =: p^{j \varepsilon}\, C_j,
$$
say. Therefore, by $\eqref{612}$, we obtain
$$
|r_m(p)| \leq  
p^{m \varepsilon}(1 + \sum_{j=2}^{N} (j-1)C_j).
$$
Hence 
$$
|\Lambda_\phi(p^m)|= |r_m(p)| \log p \leq C\,p^{m\varepsilon}\log p^m \leq C^\prime p^{m \varepsilon^\prime} .
$$
\end{pf}

From Lemma 4, we can take $\theta=\varepsilon$ for any fixed $\varepsilon >0$ in (S$5$) 
for $L_\phi(s)$ with $\phi \in C_L$. Now we define $\widetilde{S}(x)$ by 
\begin{equation}\label{617}
\widetilde{S}(x) :=\widetilde{S}_\phi(x;h) 
:=\sum_{p} \sum_{m=1}^2 \Lambda_{\phi}(p^m) h( x p^m ).
\end{equation} 

\begin{lem} Let $\phi \in C_L(N)$. 
Then 
\begin{equation}\label{618}
\aligned
\widetilde{S}(x) 
& = \sum_{p} \log p \sum_{m=1}^2 \phi(p^m) h( x p^m )
  + \sum_{p \not\in S_\phi} C(p)h(xp^2) \log p  
  + \sum_{p \in S_\phi} C(p)h(xp^2) \log p  \\
& =: \widetilde{S}_1(x) +\widetilde{S}_2(x) +\widetilde{S}_3(x)
\endaligned 
\end{equation}
say, where
\begin{equation}\label{619}
C(p) =
\begin{cases}
c_2(p) & \text{if $n_p \geq 2$}, \\
0      & \text{if $n_p=1$}.
\end{cases} 
\end{equation} 
\end{lem}

\begin{pf}
This is a direct consequence of $\eqref{611}$ and $\eqref{612}$. 
\end{pf}


\begin{lem} 
Let $\phi \in C_L(N)$ and let $\widetilde{S}_1(x)$ be as above. Then
\begin{equation}\label{620}
\widetilde{S}_1(x)
=\int_0^\infty h(xu) \phi(u) du 
-\sum_{\zeta^\ast(\rho)=0} \int_0^\infty h(xu) \phi(u) u^\rho \frac{du}{u} 
+O(x^{-\frac{1}{3}-\varepsilon})
\end{equation}
for any $\varepsilon >0$. 
\end{lem}


\begin{pf}
First we show that 
\begin{equation}\label{621}
\widetilde{S}_1(x)
= \sum_{n=1}^\infty \Lambda(n)\phi(n)h(xn)
 +O( x^{-\frac {1}{3} -\varepsilon } )
\end{equation}
for any $\varepsilon >0$. 
By a way similar to the proof of Lemma 3, we obtain 
\begin{equation}\label{622}
\sum_{p} \sum_{m=3}^\infty \phi(p^m) h( x p^m )
= O( x^{-\frac {1}{3} -\varepsilon } )
\end{equation}
by using Lemma 1 and partial summation. 
This leads to $\eqref{621}$ because 
$$
\widetilde{S}_1(x)
= \sum_{n=1}^\infty \Lambda(n)\phi(n)h(xn) 
-\sum_{p} \sum_{m=3}^\infty \phi(p^m) h( x p^m ).
$$
Applying Proposition 1 to $\zeta(s)$ and $u \mapsto h(xu) \phi(u)$ 
we have 
\begin{equation}\label{623} 
\aligned
\sum_{n=1}^\infty & \Lambda(n) h(xn) \phi(n) \\ 
=& \int_0^\infty h(xu) \phi(u)du
 -\sum_{\zeta^\ast(\rho)=0} \int h(xu) \phi(u) u^\rho \frac{du}{u} 
 +\int_0^\infty h(xu) \phi(u) \frac{du}{u} \\
& -\sum_{n=1}^\infty  \overline{\Lambda(n)} n^{-1} h(xn^{-1})\phi(n^{-1}) \\
& -(\log \pi + \,C_E)\,h(x)\phi(1)
  -W_{1/2,0}(u \mapsto h(xu)\phi(u)).
\endaligned 
\end{equation}
As for the third term and the sixth term on the right hand side, we have 
\begin{equation}\label{624}
\aligned
\int_0^\infty & h(xu) \phi(u) \frac{du}{u} 
\, - \, W_{1/2, 0}(u \mapsto h(xu)\phi(u)) \\ 
&= \int_0^\infty h(v) \phi(v/x) \frac{dv}{v}
-\int_x^\infty h(v)\phi(v/x) 
\frac{v^2}{v^2-x^2}\frac{dv}{v} \\
&= x^2 \int_0^\infty h(v) \phi(v/x) \frac{1}{v^2-x^2} \frac{dv}{v}
-\int_0^x h(v)\phi(v/x) 
\frac{v^2}{v^2-x^2}\frac{dv}{v},
\endaligned
\end{equation}
for sufficiently small $x>0$.
Because $|\phi(u)| \leq C u^2$, 
the right hand side of $\eqref{624}$ is bounded as $x \to 0$. 
The fourth term and the fifth term on the right hand side are zero 
for sufficiently small $x>0$ because the support of $h$ is compact. 
Hence we obtain
\begin{equation}\label{625} 
\aligned
\sum_{n=1}^\infty & \Lambda(n) h(xn) \phi(n) \\ 
=& \int_0^\infty h(xu) \phi(u)du
 -\sum_{\zeta^\ast(\rho)=0} \int h(xu) \phi(u) u^\rho \frac{du}{u} 
 +O(1).
\endaligned 
\end{equation}
Lemma 6 follows from $\eqref{621}$ and $\eqref{625}$. 
\end{pf}


\begin{lem} 
Let $\phi \in C_L(N)$ with $N \geq 2$ and let $\widetilde{S}_2(x)$ be as above. 
Then 
\begin{equation}\label{626}
\tilde{S}_2(x)=O(x^{-\frac{1}{2}-\varepsilon})
\end{equation}
for sufficiently small $x>0$. 
Further we obtain
\begin{equation}\label{627}
\widetilde{S}_2(x)
= \mu A_\phi \int_0^\infty h(xu^2) u^{\mu} \frac{du}{u} 
 +O(x^{ -\frac {\nu}{2} -\varepsilon})
\end{equation}
for any $\varepsilon >0$, when
\begin{equation}\label{628}
\sum_{p \leq T,\, p \not\in S_\phi} c_2(p) \log p = A_\phi T^{\mu} +O(T^{\nu +\varepsilon})
\end{equation}
holds for some constant $A_\phi$ and ${\rm Re}(\mu) > \nu >0$.
\end{lem}


\begin{pf}
Suppose that the support of $h(u) \in C_0^\infty$ is contained in $[a,b]$. 
Then, from $\psi(p^m) \ll p^{m \varepsilon}$, we have 
\begin{equation}\label{629}
\sum_{p \not\in S} h(xp^2) c_{2}(p) \log p
\ll \sum_{ \sqrt{\frac{a}{x}} \leq p \leq \sqrt{\frac{b}{x}} }p^{2\varepsilon} 
\ll x^{-\frac{1}{2}-\varepsilon}.
\end{equation}
This is the first assertion. 
By using partial summation we have
\begin{equation}\label{630}
\widetilde{S}_2(x)
= -\int_0^\infty 2xu h^\prime (xu^2)
 \left( \sum_{p \leq u,\, p \not\in S_\phi} c_2(p) \log p  \right) du.
\end{equation}
Applying the assumption of the Lemma we find that 
the right hand side of $\eqref{630}$ is equal to 
\begin{align*} 
-\int_0^\infty & 2xu h^\prime (xu^2)
 ( A_\phi u^{\mu} +O(u^{\nu+\varepsilon})) du \\
&= \mu A_\phi \int_0^\infty h(xu^2) u^{\mu} \frac{du}{u} 
+O(\int_0^\infty 2xu^2 |h^\prime(xu^2)|u^{\nu+\varepsilon}\frac{du}{u}) \\
&= \mu A_\phi \int_0^\infty h(xu^2) u^{\mu} \frac{du}{u} 
+O( x^{-\frac{\nu + \varepsilon}{2}} ).
\end{align*}
\end{pf}


\begin{lem} 
Let $\phi \in C_L(N)$ with $n \geq 2$ and let $\widetilde{S}_3(x)$ be as above. 
Then
\begin{equation}\label{631}
\widetilde{S}_3(x) = O(1).
\end{equation}
\end{lem}


\begin{pf}
By using partial summation we have
\begin{equation}\label{632}
\widetilde{S}_3(x)
= -\int_0^\infty 2xu^2 h^\prime (xu^2)
 \left( \sum_{p \leq u,\, p \in S_\phi} c_2(p) \log p  \right) \frac{du}{u}.
\end{equation}
From the proof of Lemma 4, we have $c_2(p) \ll  p^{2\varepsilon}$. 
Hence 
\begin{equation}\label{633}
\sum_{p \leq u,\, p \in S_\phi} c_2(p) \log p   \ll
\sum_{p \in S_\phi} p^{2 \varepsilon} \log p  =O(1).
\end{equation}
From $\eqref{632}$ and $\eqref{633}$ we have $\widetilde{S}_3(x) =O(1)$.
\end{pf}


\subsection{Proof of Theorem 1.}
Let $\phi \in C_L(1)$.  
From the definition of $C_L(1)$, $L_\phi(s)$ belongs to ${\cal S}$ with  $\theta=\varepsilon$.
Hence we obtain 
\begin{equation}\label{634}
S(x)
=
m_\phi \cdot \int_0^\infty h(xu) du 
-\sum_{L_\phi^\ast(\rho)=0} \int_0^\infty h(xu) u^\rho \frac{du}{u}+O(1)
\end{equation} 
by Lemma 2. 
Here we note that 
$\Lambda_{\phi}(n)=\Lambda(n)\phi(n)$ for $\phi \in C_L(1)$. 
Therefore 
\begin{equation}\label{635}
S(x) =\sum_{n=1}^\infty \Lambda(n)h(xn)\phi(n).
\end{equation}
By applying Proposition 1 to $\zeta(s)$ and $u \mapsto h(xu)\phi(u)$, we obtain
\begin{equation}\label{636}
S(x) = \int_0^\infty h(xu) \phi(u) du  
- \sum_{\zeta^\ast(\rho)=0} \int_0^\infty h(xu) \phi(u) u^\rho \frac{du}{u} + O(1),
\end{equation}
in a way similar to the proof of Lemma 2. 
Theorem 1 follows from $\eqref{634}$ and $\eqref{636}$. $\Box$ 


\subsection{Proof of Theorem 2.}
Let $\phi \in C_L(N)$ with $N \geq 2$. 
From the definition of $C_L(N)$ and Lemma 4, $L_\phi(s)$ belongs to ${\cal S}$ with  $\theta=\varepsilon$.
Hence we obtain 
\begin{equation}\label{637}
S(x)
=m_\phi \cdot \int_0^\infty h(xu) du 
-\sum_{L_\phi^\ast(\rho)=0} \int_0^\infty h(xu) u^\rho \frac{du}{u}
+O(1)
\end{equation}
by Lemma 2. On the other hand we have
\begin{equation}\label{638}
S(x) = \widetilde{S}(x)+O(x^{-\frac{1}{3}-\varepsilon})
\end{equation}
by Lemma 3. 
Together with Lemmas $5,6,7$ and $8$ we obtain 
\begin{equation}\label{639}
S(x) = 
\int_0^\infty h(xu) \phi(u) du  
- \sum_{\zeta^\ast(\rho)=0} \int_0^\infty h(xu) \phi(u) u^\rho \frac{du}{u} 
+ J_{\phi}(x,h) + O(x^{-\frac 13 -\varepsilon})
\end{equation}
where $J_{\phi}(x,h)=\tilde{S}_2(x)$. 
Theorem 2 follows from $\eqref{637}$, $\eqref{639}$ and Lemma 7. $\Box$ 

%
%
\section{ Proof of Theorem 3 and Theorem 4.}

Theorem 3 and Theorem 4 are proved by an argument quite similar to the proof of Theorem 1 and Theorem 2, 
so we only describe the outline of the proof. \\

Define the numbers $\Lambda_{\phi \otimes \psi}(n)$ by
\begin{equation}\label{701}
-\frac{L_{\phi \otimes \psi}^\prime(s)}{L_{\phi \otimes \psi}(s)}=\sum_{n=1}^\infty 
\frac{ \Lambda_{\phi \otimes \psi}(n) }{n^s}.
\end{equation}
Then we have 
\begin{align}\label{702}
\aligned
\,&\Lambda_{\phi \otimes \psi}(n)= \\
&
\begin{cases}
\left( \sum_{i=1}^M \alpha_{\phi}(p,i)^m \right) 
\left( \sum_{j=1}^N \alpha_{\psi}(p,j)^m \right) \log  p 
 &\text{if $n=p^m$ with $p \not\in S$,  $m \geq 1$,} \\ 
0&\text{if $n$ is not a power of a prime,}
\end{cases}
\endaligned
\end{align}
and the following Lemma 9 is proved similarly to Lemma 1. 


\begin{lem}
Let $(\phi, \psi) \in C_L(M) \times C_L(N)$ be a Rankin-Selberg pair. 
Then
\begin{align}\label{703}
\sum_{p \not\in S} \sum_{{l \leq m}\atop{p^m\leq T}} 
\Lambda_{\phi \otimes \psi}(p^m)
\ll_{\varepsilon,l} T^{\frac{1}{l}+\varepsilon }
\end{align}
for $l \geq 1$, and 
\begin{align}\label{704}
\sum_{p \in S} \sum_{m=1}^\infty 
\Lambda_{\phi \otimes \psi}(p^m)=O(1).
\end{align}
\end{lem}
\begin{pf}
First we note the relation
$$
b_\phi(p^m)= \frac{\Lambda_{\phi}(p^m)}{\log p^m}=\frac{1}{m} \sum_{i=1}^M \alpha_{\phi}(p,i)^m 
$$
for $p \not\in S_{\phi}$ and 
the estimate $b_\phi(n) \ll_{\varepsilon} n^{\varepsilon}$ for $\phi \in C_L$ which follows by Lemma 4.
From these we have 
\begin{equation}\label{705}
\Lambda_{\phi \otimes \psi}(p^m) =
\left( \sum_{i=1}^M \alpha_{\phi}(p,i)^m \right) 
\left( \sum_{j=1}^N \alpha_{\psi}(p,j)^m \right)
\log  p  \ll_{\varepsilon} p^{m\varepsilon}
\end{equation}
for $p \not\in S$, $(\phi, \psi) \in C_L(M) \times C_L(N)$. 
Hence $\eqref{703}$ is obtained by the same arguments as in the proof of Lemma 4.   
Because $L_{\phi \otimes \psi}(s) \in {\cal S}$,
\begin{equation}\label{706}
\Lambda_{\phi \otimes \psi}(p^m) = b_{\phi \otimes \psi}(p^m) \log p^m
\ll p^{m\theta}\log p^m \ll p^{m(\theta+\varepsilon)},
\end{equation}
where $b_{\phi \otimes \psi}(n)$ are given by 
$\log L_{\phi \otimes \psi}(s)=\sum_{n=1}^\infty b_{\phi \otimes \psi}(n)n^{-s}$.
Therefore 
\begin{equation}\label{707}
\sum_{p \in S} \sum_{m=1}^\infty \Lambda_{\phi \otimes \psi}(p^m)
 \ll \sum_{p \in S} \sum_{m=1}^\infty p^{m(\theta+\varepsilon)}
 = \sum_{p \in S} \frac{ p^{\theta+\varepsilon} }{1-p^{\theta+\varepsilon}} 
 =O(1).
\end{equation}
\end{pf}

For any fixed $h\in C_0^\infty$ and $x>0$ we consider the sum 
\begin{equation}\label{708}
S(x):=S_{\phi,\psi}(x;h):=\sum_{n=1}^\infty \Lambda_{\phi \otimes \psi}(n) h(xn).
\end{equation}
Theorem 3 or Theorem 4 is proved by computing the sum $S(x)$ in two ways. 
Applying Proposition 1 to $L_{\phi \otimes \psi}(s)$ and $u \mapsto h(xu)$, we have
\begin{equation}\label{709}
S(x)
=
m_{\phi \otimes \psi} \int_{0}^{\infty} h(xu)du
-\sum_{L_{\phi \otimes \psi}^\ast (\rho)=0} 
\int_{0}^{\infty} h(xu) u^{\rho} \frac{du}{u} +O(1) 
\end{equation}
for sufficiently small $x>0$.
By using Lemma 9 and partial summation, we obtain 
\begin{equation}\label{710}
S(x)= \sum_{p \not\in S} \sum_{m=1}^2 \Lambda_{\phi \otimes \psi}(p^m) h(xp^m) + O(x^{-\frac{1}{3}-\varepsilon}).
\end{equation} 
Take
\begin{equation}\label{711}
\widetilde{S}(x) := \sum_{p \not\in S} \sum_{m=1}^2 \Lambda_{\phi \otimes \psi}(p^m) h(xp^m),  
\end{equation} 
and divide the sum $\widetilde{S}(x)$ into two parts as 
\begin{equation}\label{712}
\aligned
\widetilde{S}(x)
&=\sum_{p \not\in S} \sum_{m=1}^2 \Lambda_\phi(p^m) h(xp^m) \psi(p^m)\log p
+ \sum_{p \not\in S} \Lambda_\phi (p^2) h(xp^2) c_{\psi,2}(p) \\
&=: \widetilde{S}_{1}(x) + \widetilde{S}_{2}(x),
\endaligned
\end{equation}
say. For $\widetilde{S}_{1}(x)$ we obtain 
$$
\widetilde{S}_{1}(x)=
\sum_{n = 1}^{\infty} \Lambda_\phi(n) h(xn) \psi(n) + O(x^{ -\frac{1}{3} -\varepsilon})
$$
by Lemma 9. 
Applying Proposition 1 to $L_\phi(s)$ and $u \mapsto h(xu)\psi(u)$, 
we obtain 
\begin{equation}\label{713}
\widetilde{S}_{1}(x)
=
m_\phi \int_{0}^{\infty} h(xu)\psi(u)du
-\sum_{L_\phi^\ast(\rho)=0} \int_0^\infty h(xu) \psi(u) u^\rho \frac{du}{u}
+O(x^{-\frac{1}{3}-\varepsilon}).
\end{equation}
For $\widetilde{S}_{2}(x)$ we find that 
\begin{equation}\label{714}
\tilde{S}_2(x)=O(x^{-\frac{1}{2}-\varepsilon})
\end{equation}
for sufficiently small $x>0$. 
In fact, when the support of $h(u) \in C_0^\infty$ is contained in the interval $[a,b]$,  
\begin{equation}\label{715}
\sum_{p \not\in S} \Lambda_\phi (p^2) h(xp^2) c_{\psi,2}(p) 
\ll \sum_{ \sqrt{\frac{a}{x}} \leq p \leq \sqrt{\frac{b}{x}} }p^{2\varepsilon} 
\ll x^{-\frac{1}{2}-\varepsilon} 
\end{equation}
by Lemma 4 and $\psi(p^m) \ll p^{m \varepsilon}$.
This implies $\eqref{714}$.

Here we note that    
\begin{equation} \label{716}
\aligned
\Lambda_\phi(p^2)
& = \left( \sum_{i=1}^M \alpha_{\phi}(p,i)^2 \right) \log p \\
& = 
\left[ 
\left( \sum_{i=1}^M \alpha_{\phi}(p,i) \right)^2 
-2 \sum_{1 \leq i<j \leq M} \alpha_{\phi}(p,i) \alpha_{\phi}(p,j)  
\right]
\log p \\
& =\Lambda_{\phi \otimes \phi}(p) + 2 c_{\phi,2}(p)\log p.
\endaligned \end{equation}
By using $\eqref{716}$ we divide $\widetilde{S}_{2}(x)$ into two parts as 
\begin{align*}
\widetilde{S}_{2}(x)
& =\sum_{p \not\in S} \Lambda_{\phi \otimes \phi}(p) h(xp^2)c_{\psi,2}(p) 
  + 2\sum_{p \not\in S} h(xp^2) c_{\phi,2}(p) c_{\psi,2}(p) \log p \\
&=:\widetilde{S}_{3}(x) + 2\widetilde{S}_{4}(x),
\end{align*}
say. 
By partial summation $\widetilde{S}_{3}(x)$ is expressed as 
\begin{equation}\label{717}
\widetilde{S}_{3}(x)
=-\int_0^\infty 
\left( \sum_{p \leq u} \Lambda_{\phi \otimes \phi}(p) c_{\psi,2}(p) \right)
(h(xu^2))^\prime du.
\end{equation}
If $\eqref{thm4:1}$ holds, we have
\begin{equation}\label{718}
\widetilde{S}_{3}(x)= \mu_1 A_1 
\int_0^\infty h(xu^2) u^{\mu_1} \frac{du}{u}
+O(x^{-\frac{\nu_1}{2}-\varepsilon})
\end{equation}
by substituting $\eqref{thm4:1}$ into $\eqref{717}$. 
Also, by using the integral expression 
\begin{equation}\label{719}
\widetilde{S}_{4}(x) = 
-\int_0^\infty \left(\sum_{p \leq u} c_{\phi,2}(p) c_{\psi,2}(p) \log p \right) ( h(xu^2) )^\prime du
\end{equation}
we obtain 
\begin{equation}\label{720} 
\widetilde{S}_{4}(x)= \mu_3 A_3 \int_0^\infty h(xu^2) u^{\mu_3} \frac{du}{u}
+O(x^{-\frac{\nu_3}{2} -\varepsilon}),
\end{equation}
if estimate $\eqref{thm4:3}$ holds.
Combining $\eqref{709}$, $\eqref{713}$, $\eqref{718}$ and $\eqref{720}$ we obtain the first half of Theorem 4. 
Also we obtain the first half of Theorem 3 by the same equations, if we replace $c_{\psi,2}(p)$ by $0$. 

To prove the latter half of Theorem 4, we divide $S_4(x)$ into four parts as  
\begin{equation}\label{721}
\aligned
\widetilde{S}(x)
=&\sum_p \sum_{m=1}^2 h(xp^m)\phi(p^m)\psi(p^m)\log p + \sum_p h(xp^2)\phi(p^2)c_{\psi,2}(p)\log p \\
 &+\sum_p h(xp^2)\psi(p^2)c_{\phi,2}(p)\log p +\sum_p h(xp^2)c_{\phi,2}(p)c_{\psi,2}(p)\log p \\
=&\sum_p \sum_{m=1}^2 h(xp^m)\phi(p^m) \psi(p^m) \log p + \sum_p h(xp^2) \phi(p)^2 c_{\psi,2}(p)\log p \\
 &+\sum_p h(xp^2) \psi(p)^2 c_{\phi,2}(p) \log p +3\sum_p h(xp^2)c_{\phi,2}(p)c_{\psi,2}(p)\log p \\
=&\sum_p \sum_{m=1}^2 h(xp^m)\phi(p^m) \psi(p^m) \log p + \sum_p h(xp^2) \Lambda_{\phi \otimes \phi}(p) c_{\psi,2}(p) \\
 &+\sum_p h(xp^2) \Lambda_{\psi \otimes \psi}(p) c_{\phi,2}(p) +3\sum_p h(xp^2)c_{\phi,2}(p)c_{\psi,2}(p)\log p \\
=&:\widetilde{S}_{5}(x)+\widetilde{S}_{6}(x)+\widetilde{S}_{7}(x)+\widetilde{S}_{8}(x),
\endaligned
\end{equation}
say. 
By a way similar to the proof of $\eqref{714}$ we obtain
\begin{equation}\label{722} 
\widetilde{S}_{l}(x)=O(x^{-\frac{1}{2}-\varepsilon}) \quad (l =6,7,8) 
\end{equation}
for sufficiently small $x>0$. Further we find that  
\begin{equation}\label{723}
\aligned 
\widetilde{S}_{6}(x) &= \mu_1 A_1 \int_0^\infty h(xu^2) u^{\mu_1} \frac{du}{u}
+O(x^{-\frac{\nu_1}{2} -\varepsilon}), \\
\widetilde{S}_{7}(x) &= \mu_2 A_2 \int_0^\infty h(xu^2) u^{\mu_2} \frac{du}{u}
+O(x^{-\frac{\nu_2}{2} -\varepsilon}) \\
\widetilde{S}_{8}(x) &= 3 \mu_3 A_3 \int_0^\infty h(xu^2) u^{\mu_3} \frac{du}{u}
+O(x^{-\frac{\nu_3}{2} -\varepsilon}) 
\endaligned 
\end{equation}
hold, if the estimates $\eqref{thm4:1}$, $\eqref{thm4:2}$ and $\eqref{thm4:3}$
hold respectively. 
We omit the process of calculations  for $ \widetilde{S}_{6}(x)$, $\widetilde{S}_{7}(x)$ and $\widetilde{S}_{6}(x)$, 
because thay are calculated in almost the same way as that for $\widetilde{S}_{2}(x)$ or $\widetilde{S}_{4}(x)$. 
For $\widetilde{S}_{5}(x)$ we have 
\begin{equation}\label{724}
\widetilde{S}_{5}(x)=\sum_{n=1}^\infty \Lambda(n)\phi(n)\psi(n)h(xn)+O(x^{-\frac{1}{3}-\varepsilon}).
\end{equation}
By applying Proposition 1 to $\zeta(s)$ and $u \mapsto h(xn)\phi(n)\psi(n)$ we have 
\begin{equation}\label{725}
\aligned
\sum_{n = 1}^\infty & \Lambda(n)\phi(n)\psi(n)h(xn)  \\
&=\int_{0}^{\infty} h(xu)\phi(u)\psi(u)du
-\sum_{\zeta^\ast(\rho)=0}
\int_0^\infty h(xu) \phi(u) \psi(u) u^{\rho} \frac{du}{u}
+O(1)
\endaligned
\end{equation}
for sufficiently small $x>0$. 
Combining $\eqref{709}$, $\eqref{723}$ and $\eqref{725}$ we obtain the latter half of Theorem 4. 
Also we obtain the latter half of Theorem 3 by the same equations, if we replace $c_{\psi,2}(p)$ by $0$. 
\hfill $\Box$

%
%
\section{Additional Topics}

\subsection{Explicit equations}
Our theorems in $\S 3$ are asymptotic results. 
We can also obtain a result which is an explicit version of our theorems in $\S3$, 
if we use an interpolation function of 
\begin{equation}\label{801}
\omega_\phi(n):=
\begin{cases}
\sum_{i=1}^{n_p} \alpha_{\phi}(p,i)^m & \text{if $n=p^m$ with $m \geq 1$,} \\
0 & \text{otherwise}.
\end{cases}
\end{equation}
However there is a possibility that such interpolation functions are not so useful for applications. 
At least it seems that 
a well-chosen interpolation function of Dirichlet coefficients is more useful than 
an interpolation function of $\eqref{801}$ for some specific purposes. 
This is one reason why we adopt asymptotic formulas as main results.  

Anyway we will establish our explicit identities. 
The key of the following results are the equations 
\begin{align}\label{802}
\Lambda_\phi(n) &= \Lambda(n)\omega_{\phi}(n), \\
\Lambda_{\phi \otimes \psi}(n) &=\Lambda_{\phi}(n)\omega_\psi(n) = \Lambda(n) \omega_\phi(n) \omega_\psi(n).
\end{align}
\begin{thm} 
Let $\phi \in C_L(M)$ and 
let $\Omega_{\phi}(u)$ be an interpolation function of $\omega_{\phi}(n)$, 
that is, $\Omega_{\phi}(n)=\omega_{\phi}(n)$ for any non-negative integer $n$. 
Then we have the following explicit identitiy 
\begin{align}
\aligned
 m_\phi  \, \widehat{h}(0) & -\sum_{L_\phi^\ast(\rho)=0}\widehat{h}(\rho) + m_\phi \, \widehat{h}(1) \\
 & -\sum_{n=1}^\infty \overline{\Lambda_\phi(n)} n^{-1} h(n^{-1})  
  -( 2 \log Q_\phi+d_\phi \,C_E)\,h(1)-\sum_{j=1}^{r_\phi} W_{\lambda_j(\phi),\mu_j(\phi)}(h), \\
=\, \widehat{ h_\Omega }(0) & -\sum_{\zeta^\ast(\rho)=0}\widehat{ h_\Omega }(\rho) +   \widehat{ h_\Omega }(1) \\
 & -\sum_{n=1}^\infty \Lambda(n) n^{-1} h_\Omega (n^{-1})  
   -(\log \pi+\,C_E)\,h_\Omega (1) - W_{\frac{1}{2},0}( h_\Omega ),
\endaligned
\end{align}
where $h_\Omega(u):=h(u)\Omega_{\phi}(u)$ and $W_{\lambda,\mu}(\cdot)$ is the functional defined in $\eqref{505}$. 
\end{thm} 
\begin{thm}
Let $(\phi,\psi) \in C_L(M) \times C_L(N)$ be a Rankin-Selberg pair.
Then we have the following explicit identities
\begin{align}
\aligned
 m_{\phi \otimes \psi} \, \widehat{h}(0) 
& -\sum_{L_{\phi \otimes \psi}^\ast(\rho)=0}
\widehat{h}(\rho) + m_{\phi \otimes \psi} \, \widehat{h}(1) \\
 & -\sum_{n=1}^\infty \overline{\Lambda_{\phi \otimes \psi}(n)} n^{-1} h(n^{-1})  
  -( 2 \log Q_{\phi \otimes \psi}+d_{\phi \otimes \psi} \,C_E)\,h(1)
  -\sum_{j=1}^{r_{\phi \otimes \psi}} W_{\lambda_j(\phi \otimes \psi),\mu_j( \phi \otimes \psi)}(h), \\
=\, m_\phi  \, \widehat{h_{\Omega_\psi}}(0) & -\sum_{L_\phi^\ast(\rho)=0}\widehat{h_{\Omega_\psi}}(\rho) 
+ m_\phi \, \widehat{h_{\Omega_\psi}}(1) \\
 & -\sum_{n=1}^\infty \overline{\Lambda_\phi(n)} n^{-1} h_{\Omega_\psi} (n^{-1})  
  -( 2 \log Q_\phi+d_\phi \,C_E)\,h_{\Omega_\psi}(1)-\sum_{j=1}^{r_\phi} W_{\lambda_j(\phi),\mu_j(\phi)}(h_{\Omega_\psi} ), \\
=\, \widehat{ h_{\Omega_\phi \Omega_\psi} }(0) & -\sum_{\zeta^\ast(\rho)=0}\widehat{ h_{\Omega_\phi \Omega_\psi} }(\rho)
 +   \widehat{ h_{\Omega_\phi \Omega_\psi} }(1) \\
 & -\sum_{n=1}^\infty \Lambda(n) n^{-1} h_{\Omega_\phi \Omega_\psi} (n^{-1})  
   -(\log \pi+\,C_E)\,h_{\Omega_\phi \Omega_\psi}(1) - W_{\frac{1}{2},0}( h_{\Omega_\phi \Omega_\psi} ),
\endaligned
\end{align}
where $h_{\Omega_\phi}(u):=h(u)\Omega_\phi(u)$ and $h_{\Omega_\phi \Omega_\psi} (u):=h(u)\Omega_\phi(u)\Omega_\psi(u)$. 
\end{thm}
Theorem 7 is obtained by calculating the sum $\sum_{n=1}^\infty  \Lambda_\phi(n) h(n)$ in two ways. 
Theorem 8 is obtained by calculating the sum $\sum_{n=1}^\infty  \Lambda_{\phi \otimes \psi}(n) h(n)$ in three ways. 
These processes are very similar to the proofs of Theorem 1 $\sim$ Theorem 4, 
therefore we omit the details of their proofs.  
We deal with one way to construct $\Omega(u)$ in the next section.  

\subsection{ One way of the construction of an interpolation function.}
In this part we give a way to construct an interpolation function by using Fourier series.
Let $a:{\Bbb N} \to {\Bbb C}$ be a function on natural numbers. 
When $a(\cdot)$ has polynomial order, 
we define the function $f_a$ by 
\begin{equation}
f_{a}(z):=\sum_{n=1}^\infty a(n) e^{2\pi i n z }.
\end{equation}
Since $a(\cdot)$ has polynomial order, 
$f_a(z)$ converges absolutely on the upper half-plane ${\rm Im}(z)>0$.  
By using $f_a$, we define $A(u)$ by 
\begin{equation}
A(u):=A(u;y,\eta):= e^{2\pi uy} \int_\eta^{\eta+1} f_a(x+iy) e^{-2\pi i u x}dx,
\end{equation}
for some fixed $y>0$ and $\eta \in {\Bbb R}$, or
\begin{equation}
A(u):=A(u;\eta):= e^{2\pi} \int_\eta^{\eta+1} f_a(x+iu^{-1}) e^{-2\pi i u x}dx,
\end{equation}
for some fixed $\eta \in {\Bbb R}$. 
From the definition, $A(u)$ satisfies $A(n)=a(n)$.

For $\phi \in C_L(n)$, we can construct the interpolation function $\Phi(u)$ of $\phi(n)$,  
since $\phi(n) \ll_{\varepsilon} n^{\varepsilon}$. 
Of course $\Phi \not\equiv \phi$ as a function on $(0,\infty)$ in general. 
Similarly we can construct the interpolation function $\Omega_{\phi}(u)$ of $\omega_\phi(n)$, 
since $\omega(n) \ll_\varepsilon n^{\theta+\varepsilon} $ for any fixed $\varepsilon >0$ from Lemma 4. 

\subsection{Symmetries of zero-sums.} 
In $\S 8.2$, we gave one way to construct an interpolation function. 
However, there is no reason that the interpolation in $\S 8.2$ is a canonical one. 
Actually, there are infinity many possibilities of interpolation functions, 
when we restrict them to the class of smooth functions. 
However the non-existence of canonical interpolation is not unfortunate. 
The existence of several different interpolation functions gives a symmetry of zero-sums.  
Let $\phi_\chi(\cdot)$, $\psi_\chi(\cdot)$ be two different interpolation functions 
of a primitive Dirichlet character $\chi$ mod $q$, 
let $\phi_f(\cdot)$, $\psi_f(\cdot)$ be two different interpolation functions 
of Fourier coefficients of $ f \in S_k(1)$ 
and let $\phi_g(\cdot)$, $\psi_g(\cdot)$ be two different interpolation functions 
of Fourier coefficients of $ g \in S_k(1)$. 
They give several ``symmetries'' of the sums $\sum_{\zeta^\ast(\rho)=0}$, $\sum_{L^\ast(\rho,\cdot)=0}$ etc.  
For example, for a suitable test function $h \in C_0^\infty$, we have 
\begin{equation*}
\aligned
\sum_{\zeta^\ast(\rho)=0} 
\int_0^\infty h(xu)\phi_\chi(u) u^\rho \frac{du}{u} 
\sim 
\sum_{\zeta^\ast(\rho)=0}
\int_0^\infty h(xu)\psi_\chi(u) u^\rho \frac{du}{u},
\endaligned
\end{equation*}
\begin{equation*}
\aligned
\sum_{\zeta^\ast(\rho)=0} 
\int_0^\infty h(xu) \phi_f(u)u^\rho \frac{du}{u} 
\sim 
\sum_{\zeta^\ast(\rho)=0}
\int_0^\infty h(xu) \psi_f(u)u^\rho \frac{du}{u},
\endaligned
\end{equation*}
and 
\begin{equation*}
\aligned
\sum_{\zeta^\ast(\rho)=0} &
 \int_0^\infty h(xu)\phi_f(u)\phi_g(u)u^\rho \frac{du}{u} 
 \sim \sum_{\zeta^\ast(\rho)=0} \int_0^\infty h(xu) \phi_f(u) \psi_g(u) u^\rho \frac{du}{u} \\ 
& \sim \sum_{\zeta^\ast(\rho)=0} \int_0^\infty h(xu) \psi_f(u) \phi_g(u) u^\rho \frac{du}{u} 
 \sim \sum_{\zeta^\ast(\rho)=0} \int_0^\infty h(xu) \psi_f(u) \psi_g(u) u^\rho \frac{du}{u},
\endaligned 
\end{equation*}
for sufficiently small  $x > 0$.
Furthermore we find that  
\begin{equation*}
\aligned
\sum_{\zeta^\ast(2\rho)=0} \int_0^\infty h(xu) u^\rho \frac{du}{u} 
\sim 
 \sum_{\zeta^\ast(\rho)=0} \int_0^\infty h(xu) \phi_\theta(u) u^\rho \frac{du}{u}
\endaligned
\end{equation*}
where $\phi_\theta$ is an interpolation function of the coefficients of 
$\theta(z)=\sum_{n=1}^\infty e^{i n^2 z}$, and  
\begin{equation*}
\aligned
\sum_{L^\ast(\rho,f)=0}
& \int_0^\infty h(xu) \phi_g(u) u^\rho \frac{du}{u}  \sim \sum_{L^\ast(\rho,f)=0} 
 \int_0^\infty h(xu) \psi_g(u) u^\rho \frac{du}{u} \\ 
& \sim \sum_{L^\ast(\rho,g)=0} 
 \int_0^\infty h(xu) \phi_f(u) u^\rho \frac{du}{u} \sim \sum_{L^\ast(\rho,g)=0} 
 \int_0^\infty h(xu) \psi_f(u) u^\rho \frac{du}{u}.
\endaligned 
\end{equation*}
These are consequences of Theorem 1 $\sim$ Theorem 4. 

\subsection{An interpretation in terms of distributions.}


In this part we give an interpretation of our results in a local-global view point. 
Throughout this part, we denote by ${\cal Z}(F)$ the set of all zeros of $F$ . 

Let $\phi \in C_L(N)$ and let $L_{\phi,p}(s)^{-1}$ be the reciprocal $p$-th Euler factor of $L_{\phi}(s)$. 
Denote by $\zeta_p(s)^{-1}$ the reciprocal $p$-th Euler factor of $\zeta(s)$. 
As explained in $\S 1$, in the local point of view, 
the relation between ${\cal Z}({L_{\phi,p}^{-1}})$ and ${\cal Z}(\zeta_p^{-1})$ 
is very simple. 
For $p \not\in S_{\phi}$, the reciprocal $p$-th Euler factor $L_{\phi,p}(s)^{-1}$ is expressed as 
\begin{equation}
L_{\phi,p}(s)^{-1} =\prod_{i=1}^N (1-\alpha_{\phi}(p,i)p^{-s})
\end{equation}
Hence, if we assume the general Ramanujan conjecture, $\alpha_{\phi}(p,i)$ can be written as 
\begin{equation}
\alpha_{\phi}(p,i)=e^{i \theta_{\phi}(p,i)}
\end{equation}
for some $ \theta_{\phi}(p,i) \in {\Bbb R}$. 
This shows that ${\cal Z}({L_{\phi,p}^{-1}})$ is 
the union of $n$-piece translation of ${\cal Z}(\zeta_p^{-1})$.  
Therefore we find that 
\begin{equation}\label{811}
\sum_{ L_{\phi,p}(\rho)^{-1}=0}  \widehat{h}(\rho) 
=\sum_{i=1}^N \sum_{\zeta_p(\rho)^{-1}=0} 
\widehat{h}\left( \rho + \sqrt{-1} \,\frac{ \theta_{\phi}(p,i)}{\log p} \right) 
\end{equation}
for any $h \in C_0^\infty$. 
Similarly we obtain 
\begin{equation}\label{812}
\sum_{ L_{\phi,\infty}(\rho)^{-1}=0}  \widehat{h}(\rho) 
=\sum_{j=1}^{r_{\phi}} \sum_{\zeta_\infty(\rho)=0} 
\widehat{h}\left( \frac{\rho + \mu_j(\phi)}{\lambda_j(\phi)} \right) 
\end{equation}
where $L_{\phi,\infty}(s)^{-1}$ is the reciprocal $\Gamma$-factor of $L_\phi(s)$.  

Next we show a relation between the global zeros and the local zeros. 
From Proposition 1, we obtain
\begin{equation}\label{813}
m_\phi  \, \widehat{h}(0) -\sum_{L_{\phi}^\ast(\rho)=0}\widehat{h}(\rho) + m_\phi \, \widehat{h}(1) 
=\sum_{p} W_\phi(h,p) + W_{\phi}(h,\infty),
\end{equation}
where
\begin{equation}\label{814}
W_\phi(h,p)=\sum_{m=1}^\infty (\Lambda_\phi(p^m) h(p^m) +\overline{\Lambda_\phi(p^m)} p^{-m}h(p^{-m}))
\end{equation}
and $W_{\phi}(h,\infty)$ is given by $W_{\lambda, \mu }(\cdot)$.
Moreover, Poisson's summation formula yields   
\begin{equation}\label{815}
W_\phi(h,p)=\sum_{L_{\phi,p}(\rho)^{-1}=0} ( \widehat{h}(\rho) + \widehat{h}(1-\overline{\rho})).
\end{equation}

Together with $\eqref{811}$, $\eqref{812}$, $\eqref{813}$ and $\eqref{815}$ we obtain
\begin{align}\aligned
\widehat{h}(1) & -\sum_{L_{\phi}(\rho)=0} \widehat{h}(\rho) +\widehat{h}(0) \\
=& \sum_{p} \sum_{ L_{\phi,p}(\rho)^{-1}=0 } \left[ \widehat{h}( \rho )
 + \widehat{h} ( 1 -\rho )\right]
 + \sum_{L_{\phi,\infty}(\rho)^{-1}=0} \widehat{h}( \rho ) \\
=& \sum_{p \not\in S_{\phi} }\sum_{i=1}^N \sum_{\zeta_p(\rho)^{-1}=0} 
\left[ \widehat{h}\left( \rho + \sqrt{-1} \,\frac{ \theta_{\phi}(p,i) }{\log p} \right) 
 + \widehat{h}\left( 1-\rho - \sqrt{-1} \,\frac{ \theta_{\phi}(p,i) }{\log p} \right) \right] \\
&+ \sum_{p \in S_{\phi}} \sum_{ L_{\phi,p}(\rho)^{-1}=0 } \left[ \widehat{h}( \rho ) + \widehat{h} ( 1 -\rho )\right] 
+ \sum_{j=1}^{r_{\phi}} \sum_{\zeta_\infty(\rho)^{-1}=0} 
\widehat{h}\left( \frac{\rho + \mu_j(\phi)}{\lambda_j(\phi)} \right).
\endaligned \end{align}
This implies that there exists some relation between ${\cal Z}(L_{\phi})$ and ${\cal Z}(\zeta)$. 

Under the Riemann hypothesis for $L_\phi(s)$ and $\zeta(s)$, 
we can show that 
\begin{equation}\label{dis:1} 
\sum_{L_{\phi}^\ast (\rho)=0} \frac{X^\rho}{\rho} \sim 
\sum_{\zeta^\ast(\rho)=0} \int_0^X \phi(u)u^{\rho}\frac{du}{u}
\quad (X \to +\infty),
\end{equation}
if $\int_0^X \phi(u) du =o(\sqrt{X})$, by a way similar to the proof of Theorem 2. 
By acting the differential operator $X\frac{d}{dX}$ on both side of $\eqref{dis:1}$, 
we obtain
\begin{equation}\label{dis:2} 
\sum_{ L_{\phi}^{\ast}(\rho)=0} X^\rho \quad 
\sim
\sum_{\zeta^\ast (\rho)=0 } \phi(X) X^\rho 
\quad (X \to +\infty).
\end{equation}
Unfortunately, two series in $\eqref{dis:2}$ do not converge for any $X$, 
hence $\eqref{dis:2}$ have no strict meaning as a function in $X$. 
However, if we define the distributions as  
$$
\int_0^\infty h(u)\sum_\rho u^{\rho}\frac{du}{u}=\sum_\rho \widehat{h}(\rho)
$$
for a test function $h \in C_0^\infty$, 
the relation $\eqref{dis:2}$ can be interpreted as a relation of two distributions. 
That is, at least in the level of distribution, 
the relation $\eqref{dis:2}$ shows that 
${\mathcal Z}(L_{\phi})$ is the ``translation'' of 
${\mathcal Z}(\zeta)$ by $\phi(\cdot)$. 
The local relations $\eqref{811}$, $\eqref{812}$ are extended to
the global relation $\eqref{dis:2}$ via the interpolation function $\phi(\cdot)$. 
We may say that 
the original RH for $\zeta(s)$ implies the RH of $L_\phi(s)$ 
in the level of distribution. 
This suggests that if the original RH is false, the RH for automorphic $L$-functions is also false.    
To show this rigorously, we need to establish the Sprindzuk type theorem as in~\cite{Sp75}. 
Althogh it is one of the most important applications of our results, 
we postpone such a study to a forthcoming work. 

%
%


\bigskip

\noindent
Graduate School of Mathematics,\\
Nagoya University,\\
Chikusa-ku, Nagoya 464-8602,\\
Japan\\
e-mail address\,:\,m99009t@@math.nagoya-u.ac.jp

\end{document}